\documentclass[a4paper,12pt,reqno]{amsart}
\usepackage[T1]{fontenc}
\usepackage[utf8]{inputenc}
\usepackage{amsmath,amssymb,amsthm,mathrsfs,mathtools}
\usepackage{lmodern,microtype,graphicx,xcolor}
\usepackage[margin=28mm]{geometry}
\usepackage[hang,flushmargin]{footmisc}
\usepackage[backref=page]{hyperref}
\usepackage{xurl}
\usepackage{todonotes}
\usepackage{comment}
\newcommand{\Hom}{\mathrm{Hom}}

\newcommand{\cube}[1]{\{0,1\}^{#1}}

\newtheorem{theorem}{Theorem}[section]
\newtheorem{lemma}[theorem]{Lemma}

\newtheorem{corollary}[theorem]{Corollary}
\newtheorem{proposition}[theorem]{Proposition}

\theoremstyle{definition}
\newtheorem{defn}[theorem]{Definition}
\newtheorem{remark}[theorem]{Remark}
\newtheorem{example}[theorem]{Example}

\newcommand{\mc}{\mathcal}

\newcommand{\mf}{\mathbf}
\newcommand{\mb}{\mathbb}

\newcommand{\wh}{\widehat}
\newcommand{\wt}{\widetilde}

\newcommand{\id}{\mathrm{id}}

\newcommand{\ab}{\mathrm{Z}}

\newcommand{\tran}{\mathrm{\Theta}}

\newcommand{\Lip}{\mathrm{Lip}}

\newcommand{\q}{\mathrm{c}}
\newcommand{\ns}{\mathrm{X}}
\newcommand{\nss}{\mathrm{Y}}
\newcommand{\co}{\circ\hspace{-0.02cm}}
\newcommand{\cu}{\mathrm{C}}

\begin{document}

\vspace*{-1cm}

\title[A proof of the Jamneshan--Tao conjecture]{
    A proof of the Jamneshan--Tao conjecture \\ [0.3em] 
    by embedding nilspaces in nilmanifolds
  }
\author{Pablo Candela}
\address{Instituto de Ciencias Matem\'aticas, Calle Nicol\'as Cabrera 13-15, Madrid 28049, Spain}
\email{pablo.candela@icmat.es}

\author{Diego Gonz\'alez-S\'anchez}
\address{Universit\'e Paris Cit\'e, Sorbonne Universit\'e, CNRS, IMJ-PRG, F-75013 Paris, France}
\email{gonzalezsanchez@imj-prg.fr}

\author{Bal\'azs Szegedy}
\address{HUN-REN Alfr\'ed R\'enyi Institute of Mathematics\\ 
Re\'altanoda utca 13-15\\
Budapest, Hungary, H-1053}
\email{szegedyb@gmail.com}

\date{}

\begin{abstract}
We prove that every compact finite-rank nilspace admits a continuous injective morphism into a toral nilspace of the same step, and hence into a connected filtered nilmanifold. Combining this embedding theorem with the general nilspace inverse theorem, we establish the Jamneshan--Tao conjecture and extend its conclusion to second-countable compact abelian groups. This yields an inverse theorem for Gowers norms on all such groups, in which the correlating functions factor through morphisms into nilmanifolds. The proof relies on substantial foundations developed in nilspace theory over the past decade, particularly the general nilspace inverse theorem and the representation of compact finite-rank nilspaces as quotients of free nilspaces.
\end{abstract}
\maketitle

\section{Introduction}
Higher-order Fourier analysis arose from arithmetic combinatorics and was inspired by developments in ergodic theory \cite{gowers,hostkra,ziegler}. Since its inception at the turn of the 21st century, much of its progress has centered on the Gowers norms and their usefulness in analyzing averages of functions over linear configurations in abelian groups. Inverse theorems for these norms, and the closely related regularity theorems, describe the algebraic structure underlying higher-order Fourier analysis. At the heart of these results are classes of structured functions that play a role analogous to that of linear characters in classical Fourier analysis. These functions consist of morphisms from abelian groups into suitable nilpotent structures, composed with Lipschitz functions on those structures. Inverse theorems in the qualitative direction (i.e.\ providing no effective bounds in the theorem) include the result of Green, Tao and Ziegler for intervals of integers (or for finite cyclic groups) \cite{GTZ}, the inverse theorems for vector spaces over finite fields of Bergelson--Tao--Ziegler \cite{BTZ} and Tao--Ziegler \cite{TZ}, and the general nilspace inverse theorem for compact abelian groups\footnote{This inverse theorem holds more broadly for inverse limits of compact coset-nilspaces of finite rank.} \cite{inverse}. There has also been much progress in the quantitative direction of providing increasingly effective bounds for various special cases of the inverse theorem, including the works \cite{GM,LSS,MannersQIinv,Mili,Peluse, Prendiville}.

Nilspace theory studies spaces equipped with cubic structures generalizing cubes on abelian groups. Introduced by Antol\'in Camarena and the third-named author \cite{ACS}, building on the parallelepiped structures of Host and Kra \cite{HKpara}, this theory has developed through the treatments in \cite{Cand:Notes1,Cand:Notes2,GMV1,GMV2,GMV3}, and it provides a useful approach to understanding the Gowers norms in general \cite{couplings} (rather than in special classes of abelian groups). This approach has led to general inverse and regularity theorems for compact abelian groups, in which the underlying structures are compact finite-rank (\textsc{cfr}) nilspaces \cite[Theorem 1.6]{inverse}. Under additional hypotheses, these underlying nilspaces can be chosen to be \emph{toral}; this means that their structure groups are all tori, and implies that these nilspaces are isomorphic to connected filtered nilmanifolds (more precisely, filtered nilmanifolds with connected filtration groups) \cite[Theorem 1.7]{inverse}. 

A central conjecture in higher-order Fourier analysis, formulated by Jamneshan and Tao in
\cite[Conjecture 1.11]{JT}, proposes that filtered nilmanifolds suffice as underlying structures in an inverse theorem for arbitrary \emph{finite} abelian groups (in other words, that in the special case of \cite[Theorem 1.6]{inverse} for finite abelian groups, one can replace the underlying \textsc{cfr} nilspaces by the more specific class of filtered nilmanifolds). This conjecture has motivated various works (for instance \cite{CGSSrank,CGSS-via,CGSS-bexp,JST-Abram,JST-bexp,JST-disco}) confirming several special cases or making progress on closely related questions.

The main result of this paper confirms the Jamneshan--Tao conjecture, and extends its conclusion to second-countable compact (Hausdorff) abelian groups.\footnote{In fact, because of the more general aspect of the nilspace inverse theorem mentioned in the previous footnote, the extension of the confirmed Jamneshan-Tao conjecture holds for functions on any inverse limit of \textsc{cfr} coset nilspaces. In particular, the resulting inverse theorem for functions on nilmanifolds involves correlating functions that are also defined through nilmanifolds.}

To obtain the conjecture's desired conclusion, formulated  in terms of filtered nilmanifolds and polynomial maps, our proof  proceeds through the broader and  more stable category of topological nilspaces and their continuous morphisms. The results and tools that have been developed in nilspace theory enable a useful analysis via key constructions \emph{within} this category, including factors, extensions, fiber-products, inverse limits (see \cite{Cand:Notes1,Cand:Notes2,CGSS-p-hom, GMV1,GMV2,GMV3}), and also the more recent representation of \textsc{cfr} nilspaces as quotients of free nilspaces \cite{CGSS-doucos}. Using these tools, we obtain the following theorem, which is the central ingredient in our proof of the conjecture.

\begin{theorem}[Toral embedding]\label{thm:embedding-intro}
For any natural number $k$,  every $k$-step \textsc{cfr} nilspace admits a continuous injective nilspace morphism into a $k$-step toral  nilspace.
\end{theorem}
\noindent In this paper, ``embedding'' means a continuous injective nilspace morphism. One can also strengthen this notion and Theorem \ref{thm:embedding-intro}; see Remark \ref{rem:induced} and Theorem \ref{thm:strongembedding}.

The inverse theorem  \cite[Theorem 1.6]{inverse} provides a morphism from the initial abelian group into a \textsc{cfr} nilspace, and our embedding theorem then places this nilspace inside a filtered nilmanifold with connected filtration groups. This completes the structural task of obtaining an underlying nilspace of the desired form. Then, a short approximation argument completes the analytic task of obtaining a suitable Lipschitz function for the correlation, thus finishing the proof of the Jamneshan–Tao conjecture. We now state the resulting inverse theorem precisely.

For a compact abelian group $\ab$, we denote by $\mu_{\ab}$ its normalized Haar measure (for finite groups, this is normalized counting measure), and we write $\mb{E}_{x\in \ab}f(x)$ for the integral $\int_{\ab} f\,d\mu_{\ab}$. For $d\geq 2$ and any 
bounded Borel measurable function $f:\ab\to\mb{C}$, the Gowers $U^d$-norm of $f$ is
\[
 \|f\|_{U^d(\ab)}
 :=\Big(
 \mb{E}_{x,z_1,\ldots,z_d\in \ab}
 \prod_{v\in\{0,1\}^d}
 \mc{C}^{|v|}
 f(x+v_1z_1+\cdots+v_dz_d)\Big)^{1/2^d},
\]
where $\mc{C}$ denotes complex conjugation and
$|v|=v_1+\cdots+v_d$.
A function is $1$-bounded if its modulus is at most $1$.
For a function $F$ on a metric space $M$ with metric $d$, we denote by $\|F\|_{\Lip}$ the norm $\|F\|_\infty+\sup_{x,y\in M, x\ne y}
|F(x)-F(y)| /d(x,y)$.

\begin{theorem}[Inverse theorem for compact abelian groups]\label{thm:JT}
Let $k\geq 1$ and $0<\delta\leq 1$. There exist a finite family $\mc{N}_{k,\delta}$ of filtered nilmanifolds of degree at most $k$ with connected filtration groups, each equipped with a smooth Riemannian metric, and constants $L=L(k,\delta)$ and $\varepsilon=\varepsilon(k,\delta)>0$, with the following property. For every second-countable compact Hausdorff abelian group $\ab$ and every $1$-bounded Borel measurable function $f:\ab\to\mb{C}$ with $\|f\|_{U^{k+1}(\ab)}\geq\delta$, there are $G/\Gamma\in\mc{N}_{k,\delta}$, a continuous nilspace morphism
$\phi:\mc{D}_1(\ab)\to G/\Gamma$, and a $1$-bounded Lipschitz function $F:G/\Gamma\to\mb{C}$ with $\|F\|_{\Lip}\leq L$ such that 
$\left|\mb{E}_{x\in \ab}f(x)\overline{F(\phi(x))}\right|
 \geq\varepsilon$. 
\end{theorem}

\noindent The special case of this theorem for finite abelian groups $\ab$ implies \cite[Conjecture 1.11]{JT}, with the additional conclusion that the filtration groups are connected (note that in
\cite[Definition 1.7]{JT}, a nilspace morphism from an abelian group $\ab$ to a filtered nilmanifold is called a \emph{polynomial map}).

The proof of Theorem~\ref{thm:embedding-intro} uses the free nilspace presentation theorem 
\cite[Theorem 6.2]{CGSS-doucos}, which represents every
$k$-step \textsc{cfr} nilspace $\ns$ as a quotient $F/\Gamma$
of a free nilspace $F$ by a translation lattice $\Gamma$. The definitions
and background for these notions are recalled in
Section~\ref{sec:preliminaries}. A free nilspace has integer
and real coordinates, each assigned a degree that specifies
its cube structure. The familiar expression
$\mb{T}^d=\mb{R}^d/\mb{Z}^d$ illustrates the
1-step case of this presentation theorem. In higher step, the transformations in
$\Gamma$ may add polynomial expressions in lower-degree
coordinates to higher-degree coordinates. The defining conditions for $\Gamma$ to be a translation lattice are the algebraic property called \emph{fiber-transitivity} and two topological properties that can be viewed as nilspace-theoretic analogues of discreteness and co-compactness (called \emph{fiber-discreteness} and \emph{fiber-cocompactness}).

An important additional algebraic property that such a lattice can have is
\emph{pureness}, introduced in
\cite[Definition 5.40]{CGSS-doucos}. The relevance of this property to our strategy lies in the fact that a \textsc{cfr} nilspace
with a pure free presentation can be embedded into a toral nilspace \cite[Proposition 5.47]{CGSS-doucos}.
The embedding works by replacing the integer coordinates of the
free nilspace by real coordinates, and extending the polynomial
formulas defining the lattice action. For a pure
presentation, this gives a toral quotient into which the
original quotient embeds.

The known presentation theorem mentioned above (\cite[Theorem 6.2]{CGSS-doucos}) does not ensure pureness of $\Gamma$.
Our main structural task, carried out in
Section~\ref{sec:pure-presentations},  is therefore to strengthen this theorem by
constructing a new free nilspace $\mc F$ and a pure
translation lattice $\Gamma'$ such that
$\ns\cong\mc F/\Gamma'$. We achieve this in Theorem~\ref{thm:pure-presentation}. We prove this by induction
on the step. The $k$-step \textsc{cfr} nilspace $\ns$ is obtained from its
$(k-1)$-step factor $\ns_{k-1}$ by an extension whose fibers are
principal homogeneous spaces of a compact abelian Lie
group $\ab_k$. Assuming a pure presentation $\ns_{k-1}\cong F/\Gamma$,
we pull these fibers back to $F$ (using a standard nilspace fiber-product). Over this free nilspace $F$, the extension splits into a product nilspace $F\times \ab_k$, by \cite[Theorem 4.1]{CGSS-doucos}. We then
need to lift the action of $\Gamma$ when the 
group $\ab_k$ is replaced by a covering group of
the form $\wt{\ab}=\mb{Z}^a\times\mb{R}^b$,
while preserving all relations in the acting lattice $\Gamma$.

The difficulty is that lifting individual transformations to the covering group need not preserve the relations between them. We address this by first expressing the action in terms of functions on $\Gamma$. In the product description above, a point consists of a base point in $F$ and a fiber coordinate in $\ab_k$. We keep the base point and replace the fiber coordinate by a function that records the coordinate's value after each transformation from $\Gamma$. In this representation, the action on the function coordinate is simply a shift of its argument, and hence respects all group relations automatically. We then lift these functions to the covering group $\wt{\ab}$. Different lifts differ by functions taking values in the kernel of the covering map. By allowing suitable polynomial corrections of this kind, we obtain a family of lifted functions preserved by the action. These corrections are parametrized by finitely many integers, which become the auxiliary coordinates of the new free nilspace \(\mc F\).

In subsections \ref{subsec:invaset} and \ref{subsec:freestruct} we construct this free nilspace
and the action of $\Gamma'$ by changes of argument and
addition of correction functions. In Subsection \ref{subsec:mainpfs} we establish the translation-lattice conditions and confirm pureness. Then in Subsection \ref{subsec:endpurepres}
we prove that the quotient is isomorphic to $\ns$. Thus the
extra coordinates permit a pure action without changing
the nilspace represented by the quotient.
Proposition~5.47 of \cite{CGSS-doucos} then gives the toral embedding in Subsection \ref{subsec:endembed}. A 3-step example illustrating the proof's mechanism is presented in Subsection \ref{subsec:3step}.

Finally, Section~5 deduces Theorem~\ref{thm:JT} from these
embeddings. On each embedded image, the original
Lipschitz function becomes a continuous function. We
approximate it there by a $1$-bounded Lipschitz function
on the ambient nilmanifold, with sufficiently small uniform
error to preserve positive correlation. The inverse theorem
of \cite{inverse} provides finitely many possible
\textsc{cfr} targets and a common Lipschitz bound, so this
approximation can be made with bounds uniform in the
initial group and function.

We now turn to the
preliminaries needed for the proof.

\section{Preliminaries}\label{sec:preliminaries}

This section gathers some basic notions needed in the rest of the paper. These being mostly background on nilspace theory and related tools, this section can be skippped by any reader familiar with this background.

All topological spaces in this paper are Hausdorff and second-countable. We denote the set of integers $\{1,\ldots,n\}$ by  $[n]$ and write $\{0,1\}^n$ for the discrete $n$-cube.
\subsection{Nilspaces and filtered groups}\hfill\\
A nilspace $\ns$ is a nonempty set with specified sets $\cu^n(\ns)\subset \ns^{\{0,1\}^n}$ of cubes, satisfying three axioms \cite[Definition 1.2.1]{Cand:Notes1}. The first one is the \emph{composition axiom}, according to which, given a cube $\q\in \cu^n(\ns)$ and a \emph{discrete-cube morphism} $\phi:\{0,1\}^m\to \{0,1\}^n$ (that is, a map $\{0,1\}^m\to \{0,1\}^n$ that extends to an affine group homomorphism $\mb{Z}^m\to\mb{Z}^n$; see \cite[Definition 1.1.1]{Cand:Notes1}), we have $\q\co\phi\in \cu^m(\ns)$. The second axiom, called \emph{ergodicity}, states simply that $\cu^1(\ns)=\ns^{\{0,1\}}$. The last axiom is called \emph{corner completion}. An \emph{$n$-corner} is a configuration missing the vertex $1^n:=(1,\ldots,1)$, i.e.\ a map $\{0,1\}^n\setminus \{1^n\}\to\ns$, whose restrictions to the $(n-1)$-dimensional \emph{faces} (\cite[Definition 1.1.4]{Cand:Notes1}) of $\{0,1\}^n$ not containing $1^n$ are cubes in $\cu^{n-1}(\ns)$. The corner completion axiom requires that for every such $n$-corner $\q'$ (for every $n$), there is an cube $\q\in \cu^n(\ns)$ such that $\q(v)=\q'(v)$ for all $v\in \{0,1\}^n\setminus\{1^n\}$.

A nilspace $\ns$ is said to be a \emph{$k$-step} nilspace if every $(k+1)$-corner has a unique completion.

A \emph{compact nilspace} is a nilspace with a compact topology for which each $\cu^n(\ns)$ is a closed subset of $\ns^{\{0,1\}^n}$ with the product topology \cite[Definition 1.0.2]{Cand:Notes2}. 

An \textsc{lch} \emph{nilspace} (short for locally compact Hausdorff) is a nilspace with a locally-compact (Hausdorff) topology for which each $\cu^n(\ns)$ is closed and the restriction map $\cu^n(\ns)\to\{n\text{-corners}\}$ is open \cite[Definition 2.1]{CGSS-doucos}. 

A \emph{morphism} is a map $\varphi:\ns\to \nss$ between nilspaces that preserves cubes. A \emph{fibration} (or \emph{totally surjective morphism}) \cite[Definition 3.3.1]{Cand:Notes1} is a morphism $\varphi$ such that for any $n$-corner $\q$ in $\ns$ and any $\wt{\q}\in\cu^n(\nss)$ with $\varphi\co \q(v) = \wt{\q}(v)$ for all $v\neq 1^n$, there exists $\q'\in \cu^n(\ns)$ such that $\varphi\co\q'=\wt{\q}$ and $\q(v)=\q'(v)$ for all $v\not=1^n$. Furthermore, every morphism between \textsc{lch} or compact nilspaces is required to be continuous.

Given a $k$-step nilspace $\ns$, for each integer $i\in [0,k]$ we denote by $\ns_i$ the \emph{canonical $i$-step factor} of $\ns$ \cite[Lemma 3.2.10]{Cand:Notes1}, denoting by $\pi_i$ the corresponding factor map $\ns\to\ns_i$, with $\ns_0$ being a singleton. We denote by $\ab_i(\ns)$ the $i$-th \emph{structure group} of $\ns$ \cite[Definition 3.2.18 and Theorem 3.2.19]{Cand:Notes1}. Each factor $\ns_i$ is a \emph{degree-$i$ abelian extension} of $\ns_{i-1}$: the fibers of the map $\ns_i\to\ns_{i-1}$ are torsors (principal homogeneous spaces) for $\ab_i(\ns)$, and the differences between cubes above the same base cube are cubes in the degree-$i$ abelian nilspace $\mc{D}_i(\ab_i(\ns))$, whose definition we recall next. 

For an abelian group $\ab$ and $i\ge0$, the \emph{degree-$i$ structure} on $\ab$ is the nilspace denoted by $\mc D_i(\ab)$, in which for each $n\ge 0$ the $n$-cubes are the maps
$\q:\cube{n}\to\ab$ satisfying
\begin{equation}\label{eq:abelian-cubes}
 \sum_{v\in F}(-1)^{|v|}\q(v)=0
 \quad\text{for every $(i+1)$-dimensional face $F\subseteq\cube{n}$}
\end{equation}
where $|v|:=v_1+\cdots+v_n$.
For $n\le i$, this condition is vacuous. We also denote by $\mc D_0(\ab)$\footnote{Note that, while $\mc{D}_0(\ab)$ is not a nilspace (as it does not satisfy the ergodicity axiom), we can extend the definition of morphism to this space. Namely, we state that a morphism $\varphi:\ns\to \mc{D}_0(\ab)$ is a map that sends cubes to cubes (and it is continuous, if the spaces $\ns$ and $\ab$ are endowed with a topology). Hence $\hom(\ns,\mc{D}_0(\ab))$ is the set of constant $\ab$-valued maps on $\ns$.} the non-ergodic nilspace whose cubes are the constant maps, and let $\mc D_i(\ab)=\{0\}$ for $i<0$. 

If $\ns$ is in addition a compact (resp.\ \textsc{lch}) nilspace, we say that it is \emph{compact finite-rank}, or \textsc{cfr} for short (resp.\ \emph{Lie fibered}), if all its structure groups are compact abelian Lie groups (resp.\ compactly generated abelian Lie groups), see \cite{Cand:Notes2,CGSS-doucos}. If all structure groups of a \textsc{cfr} nilspace are connected (and hence are tori), we say that it is a \emph{toral nilspace} \cite[\S 2.9.1]{Cand:Notes2}.

We will also use the concept of a nilspace \emph{extension}, with definitions and conventions of \cite[\S 2.4 and 3.3]{Cand:Notes1}.

\medskip

A filtration on a group $G$ is a sequence of subgroups $G=G_0=G_1\ge G_2\ge\cdots$ (often denoted by $G_\bullet$) satisfying the commutator inclusions $[G_i,G_j]\subseteq G_{i+j}$ for every $i,j\geq 0$. We say that $G_\bullet$ has \emph{degree is at most $k$} if $G_{k+1}=\{1\}$. The \emph{group nilspace} associated with $G_\bullet$ is $G$ with the set of $n$-cubes being the \emph{Host--Kra cube group} $\cu^n(G_\bullet)$, i.e.\ the subgroup of $G^{\{0,1\}^n}$ generated by the maps $\{0,1\}^n\to G$ that equal $h\in G_s$ on a codimension-$s$ face and equal $1$ elsewhere \cite[Definition 2.2.3]{Cand:Notes1}. We regard filtered groups as group nilspaces when discussing morphisms. Their canonical factors and structure groups are $G/G_{i+1}$ and $G_i/G_{i+1}$, respectively, see \cite[\S 2.2]{Cand:Notes1}.

A \emph{filtered nilmanifold} of degree $d$ is a quotient space $G/\Gamma$, where $G$ is a nilpotent Lie group with a filtration $G_\bullet=(G_i)_{i\ge 0}$ of degree $d$ consisting of closed subgroups, and where $\Gamma$ is a discrete subgroup of $G$ such that $\Gamma\cap G_i$ is cocompact in $G_i$ for every $i$ \cite[\S 1.1]{Cand:Notes2}. The $n$-cubes on $G/\Gamma$ are the compositions of Host--Kra cubes in $\cu^n(G_\bullet)$ with the quotient map $G\to G/\Gamma$. By \cite[Theorem 2.9.17]{Cand:Notes2}, toral nilspaces are  filtered nilmanifolds with connected filtration groups. 

A $k$-step \emph{free nilspace} \cite[Definition 3.1]{CGSS-doucos} is a \emph{product nilspace} \cite[Definition 3.1.2]{Cand:Notes1} of the form
\begin{equation}\label{eq:free}
F=\prod_{i=1}^k\mc D_i(\ab_i),\text{ where }
 \ab_i=\mb{Z}^{a_i}\times\mb{R}^{b_i},
\end{equation}
with the product topology. 

\subsection{Translations, pureness and period lattices}\hfill\\
Let $\ns$ be a $k$-step nilspace. A translation of height $i$ on $\ns$ is a bijection $\ns\to\ns$ that preserves cubes when applied to the cube's values on any face of codimension $i$ \cite[\S 3.2.4]{Cand:Notes1} (when the nilspace has a topology, its translations are required to be homeomorphisms). Translations of height 1 are simply called \emph{translations}, and these form a group, denoted $\tran(\ns)=\tran_1(\ns)$. For each $i\in\mb{N}$, the translations of height $i$ on $\ns$ form a subgroup $\tran_i(\ns)\leq \tran(\ns)$, and the sequence $(\tran_i(\ns))_{i\ge 0}$ (with $\tran_0(\ns)=\tran(\ns)$) is a filtration \cite[Lemma 3.2.35]{Cand:Notes1}. 

The action of the translation group on the canonical factors $\ns_i$ induces homomorphisms $\eta_j:\tran(\ns)\longrightarrow\tran(\ns_j)$ defined by the property $\pi_j\co\alpha=\eta_j(\alpha)\co\pi_j$ for every $\alpha\in \tran(\ns)$ (in particular $\eta_0$ maps all of $\tran(\ns)$ to the identity). For a subgroup $\Gamma\le\tran(\ns)$ and $i\ge 0$, we shall always denote by $\Gamma_i$ the subgroup $\Gamma\cap\tran_i(\ns)$. Then $(\Gamma_i)_{i\geq 0}$ is a filtration on $\Gamma$.

Note that $
\tran(\ns)=\ker(\eta_0)\geq \ker(\eta_1) \geq \cdots \geq \ker(\eta_{k-1})\geq \ker(\eta_k)=\{\id\},$ as $\eta_k=\id$.
Thus we have the partition
\begin{equation}\label{eq:transparti}
\tran(\ns)=\{\id\}\cup\bigsqcup_{i\in [k]} \ker(\eta_{i-1})\setminus \ker(\eta_i).
\end{equation}
We also have in general $\ker(\eta_{i-1})\supset \tran_i(\ns)$ and this inclusion can be strict. The useful property called \emph{pureness}, introduced in \cite[Definition 5.40]{CGSS-doucos}, ensures that these inclusions become equalities upon intersecting with $\Gamma$. 
\begin{defn}[Pure group of translations]\label{def:pureness}
A subgroup $\Gamma\le\tran(\ns)$ is pure if
\[
 \Gamma\cap\ker\eta_{i-1}=\Gamma_i
\text{ for }1\le i\le k.
\]
\end{defn}
For a free nilspace $F$, a height-$i$ translation induces a constant displacement on the structure group $\ab_i(F)$ (which equals the group $\ab_i$ from the expression \eqref{eq:free}). This follows from the description of translations in \cite[Theorem 3.15]{CGSS-doucos}. Therefore, a non-identity element of a pure subgroup of $\tran(F)$ translates by a nonzero constant displacement in the first component of $F$ (i.e. the component $\ab_i$ in \eqref{eq:free} with least $i$) on which it acts nontrivially. In particular, a pure group acts freely. 

A subgroup $\Gamma\le\tran(\ns)$ is \emph{fiber-transitive} if, whenever $y=\gamma x$ and $\pi_{i-1}(x)=\pi_{i-1}(y)$, there exists $\gamma'\in\Gamma_i$ with $y=\gamma'x$. Pure groups on free nilspaces are fiber-transitive: if $\gamma\notin\Gamma_i$, its first nonzero constant displacement occurs below level $i$, so $y=\gamma x$ must then imply that $\pi_{i-1}(y)\neq \pi_{i-1}(x)$. Conversely, fiber-transitive groups which act freely are pure, see \cite[Proposition 5.42]{CGSS-doucos}.

\begin{defn}\label{def:translation-lattice}(Translation lattice \cite[Definition 2.13]{CGSS-via})
Let $F=\prod_{i=1}^k\mc{D}_i(\ab_i)$ be a free nilspace. A fiber-transitive group $\Gamma\le\tran(F)$ is a \emph{translation lattice} if, for each $i$, the subgroup $\eta_i(\Gamma)\cap\tran_i(F_i)$ is discrete and cocompact in $\tran_i(F_i)\cong \ab_i$.
\end{defn}
When the group $\Gamma\le \tran(F)$ is pure, we have
\begin{equation}\label{eq:layer-periods}
\eta_i(\Gamma)\cap\tran_i(F_i)=\eta_i(\Gamma_i),
\end{equation}
Indeed, the forward inclusion holds because an element of $\Gamma$ whose induced action on $F_i$ has height $i$ fixes $F_{i-1}$, and hence belongs to $\Gamma_i$ by pureness (the backward inclusion is clear). We shall often denote by $L_i$ the group $\eta_i(\Gamma_i)$ viewed as a subgroup of $\ab_i$, where an element of $\tran_i(F_i)$ is identified with the element of $\ab_i$ that it adds as a constant shift.

For a translation lattice $\Gamma\le \tran(F)$, the orbit quotient $F/\Gamma$ is a \textsc{cfr} nilspace $\ns$ with the quotient topology and the cubes projected from $F$, by \cite[Corollary 5.37]{CGSS-doucos}. The projection $F\to F/\Gamma$ is a fibration and the structure groups of $F/\Gamma$ are $\ab_i(F)/(\eta_i(\Gamma)\cap\tran_i(F_i))$, see \cite[Lemma 5.24, Theorem
5.35, and Corollary 5.37]{CGSS-doucos}.

We will use the following elementary fact about lattices.
\begin{lemma}\label{lem:lattices}
Let $a,b\geq 0$, and let $L$ be a discrete cocompact
subgroup of $\mb{Z}^a\times\mb{R}^b$. Let $H=\operatorname{pr}_{\mb{Z}^a}(L)$ and  $L_{\mb{R}}
=\{r\in\mb{R}^b:(0,r)\in L\}$. 
Then the following holds:\setlength{\leftmargini}{0.7cm}
\begin{enumerate}
\item $H$ has finite index in $\mb{Z}^a$.
\item $L_{\mb{R}}$ is a discrete cocompact subgroup of
$\mb R^b$.
\item $L$ is free abelian of rank $a+b$. Under the 
inclusion $\mb{Z}^a\times\mb{R}^b\subseteq
\mb{R}^{a+b}$, it is a full lattice in $\mb{R}^{a+b}$, i.e.
it is generated over $\mb{Z}$ by a real vector-space
basis of $\mb{R}^{a+b}$.
\end{enumerate}
\end{lemma}

\begin{proof}
Let $\ab=\mb Z^a\times\mb{R}^b$. To prove claim $(i)$, we let $\Phi:\ab/L\to\mb Z^a/H$ be $\Phi\bigl((z,r)+L\bigr)=z+H$. This map is well-defined, as changing $(z,r)$ by an element of $L$ changes $z$ by an element of $H$ and it is a continuous surjective homomorphism. Its target is discrete, so compactness of $\ab/L$ implies that $\mb{Z}^a/H$ is finite. This proves claim $(i)$.

To see claim $(ii)$, consider the map $\iota:\mb R^b/L_{\mb R}\to \ab/L$ given by $\iota(r+L_{\mb{R}})=(0,r)+L$. 
This map is well-defined, continuous, and injective. The image of $\iota$ is $\ker\Phi$. Indeed, the forward inclusion is clear, and to see the opposite inclusion note that if
$(z,r)+L\in\ker\Phi$, then $z\in H$, so by definition of $H$ there exists $(z,r_0)\in L$, and subtraction gives
$(z,r)+L=(0,r-r_0)+L$.

The map $\iota$ is also a homeomorphism onto
$\ker\Phi$. To check this, let
$q:\ab\to \ab/L$ be the quotient map. For any open
$U\subseteq\mb{R}^b$, the set $\{0\}\times U$ is open
in $\ab$, since $\mb{Z}^a$ is discrete. Quotient maps
of topological groups are open, so
$q(\{0\}\times U)=\iota(U+L_{\mb{R}})$ is open in
$\ab/L$. Thus $\iota$ is open onto its image.

The kernel of $\Phi$ is closed in the compact group $\ab/L$,
and hence compact. Therefore
$\mb{R}^b/L_{\mb{R}}$ is compact.
The subgroup $L_{\mb{R}}$ is discrete because $L$
is discrete. This proves claim $(ii)$.

Finally, to prove claim $(iii)$, we first note that by standard arguments the lattice $L_{\mb{R}}$ has a $\mb{Z}$-basis consisting of $b$ vectors $r_1,\ldots,r_b$ that form a real basis of $\mb{R}^b$. Since $H$ has finite index in $\mb{Z}^a$, it has
a $\mb{Z}$-basis $h_1,\ldots,h_a$, which is also
a real basis of $\mb{R}^a$. Choose lifts
$\ell_i=(h_i,t_i)\in L$. We claim that $\ell_1,\ldots,\ell_a,\ (0,r_1),\ldots,(0,r_b)$ form a $\mb{Z}$-basis of $L$.

To prove generation, take $(z,r)\in L$. Express
$z$ as an integer combination of the $h_i$ and
subtract the corresponding combination of the
$\ell_i$. The remaining element has integer
coordinate zero, so its real coordinate belongs
to $L_{\mb{R}}$ and is an integer combination
of the $r_j$. Independence is easy to see.

The same argument, with real coefficients, shows
that these $a+b$ vectors form a real basis of
$\mb{R}^{a+b}$. Their integer span is therefore
discrete, and the compact parallelepiped obtained
by allowing every coefficient to lie in $[0,1]$
meets every coset. Hence $L$ is a full lattice
in $\mb{R}^{a+b}$.
\end{proof}

\subsection{Two properties of translation actions}\hfill\\
The following observations will replace calculations in group coordinates.

\begin{lemma}\label{lem:action}
Let $G_\bullet$ be a filtered group acting on a nilspace $\nss$, and suppose each element of $G_s$ acts by a translation of height at least $s$. Then the action\footnote{If $G$ and $\nss$ have a topology, we assume this action to be continuous.} $G\times \nss\to \nss$ is a morphism.
\end{lemma}
\begin{proof}
Let $\q\in \cu^n(G)$ and $\q^*\in \cu^n(\nss)$. Express $\q$ as a product of the face generators defining $\cu^n(G)$ (see \cite[Lemma 2.2.5]{Cand:Notes1}). Applying such a generator pointwise to $\q^*$ is precisely an allowed face translation. Applying the factors successively proves that $v\mapsto \q(v)\q^*(v)$ is a cube.
\end{proof}

Observe that the space of morphisms $\Hom(\mc{D}_1(\mb{R}),\mc D_d(\mb{R}))$ is the set of classical polynomial functions of degree at most $d$ from $\mb{R}$ to $\mb{R}$ (this follows from \cite[Lemma 3.6]{CGSS-doucos}). In our upcoming constructions for the proof of Theorem \ref{thm:JT}, we will need to work with the following ample and natural generalization of this notion of polynomials where the domain is allowed to be any nilspace, and the codomain any abelian group. 
\begin{defn}[Polynomials on nilspaces]\label{def:nilpolys}
For a nilspace $S$, an abelian group $\ab$, and an integer $d$, let
\begin{equation}\label{eq:polynomials}
 \mc{P}_d(S,\ab):=\Hom(S,\mc D_d(\ab)).
\end{equation}
\end{defn}
This is an abelian group under pointwise addition. Note that, while this could seem to be just a special case of the notion of morphisms between nilspaces, this specific case involves the \emph{degree} parameter $d$, which induces useful filtrations on this abelian group.

\begin{remark}\label{rem:polyfiltra}
For any $d$, the standard filtration that we use on $\mc P_d(S,\ab)$ is $(\mc P_d(S,\ab))_i:=\mc P_d(S,\ab)$ for $i=0,1$ and $(\mc P_d(S,\ab))_i:=\mc P_{d-i+1}(S,\ab)$ for $i\in[2,d+1]$ (this is a filtration of degree $d+1$). Note that, as group nilspaces, we have that $\mc P_d(S,\ab)$ is a subnilspace (see Remark \ref{rem:induced}) of $\mc P_{d'}(S,\ab)$ whenever $d<d'$. However $\mc P_d(S,\ab)$ is not an \emph{induced} subnilspace of $\mc P_{d'}(S,\ab)$ (see Remark \ref{rem:induced}), since $\cu^n(\mc P_{d'}(S,\ab))\cap \mc P_d(S,\ab)^{\{0,1\}^n}$ in general can \emph{strictly} include $\cu^n(\mc P_d(S,\ab))$ (essentially because the filtration on $P_{d'}(S,\ab)$ has greater degree). 
\end{remark}

As we noted before, if $S$ and $\ab$ have topologies, we define \eqref{eq:polynomials} as the subgroup of continuous morphisms, with the compact-open topology. Note that if $S$ is discrete, every map in $\mc{P}_d(S,\ab)$ is continuous (and this will be the case in this paper). For $d<0$ we take $\mc{P}_d(S,\ab)$ to be the zero group, and $\mc{P}_0(S,\ab)$ consists of constant functions.

\begin{lemma}[Difference under a translation]\label{lem:difference}
If $\alpha\in\tran_s(S)$ and $q\in\mc P_d(S,\ab)$, then $q\co\alpha-q\in\mc P_{d-s}(S,\ab)$.
\end{lemma}
\begin{proof}
For a cube $\q\in \cu^n(S)$, consider the arrow $\langle \q,\alpha\co\q\rangle_s$ (see \cite[Definition 2.2.18]{Cand:Notes1}). By definition of $\alpha$ being in $\tran_s(S)$ and \cite[Lemma 3.2.32]{Cand:Notes1}, we have $\langle \q,\alpha\co\q\rangle_s\in \cu^{n+s}(S)$. Hence $q\co \langle \q,\alpha\co\q\rangle_s=\langle q\co \q,q\co \alpha\co\q\rangle_s$ is in $\cu^{n+s}(\mc{D}_d(\ab))$. By \cite[Lemma 3.3.37]{Cand:Notes1} this implies that $(q\co\alpha-q)\co \q\in\cu^n(\mc{D}_{d-s}(\ab))$. If $s>d$, it is zero. Since this holds for every $\q$, we conclude that $q\co\alpha-q\in\mc P_{d-s}(S,\ab)$ as claimed.
\end{proof}
We record the following useful consequence of Lemma \ref{lem:difference}. For a filtered group $G_\bullet$, right multiplication by $h\in G_s$ is a translation of height at least $s$ on its group nilspace. Then Lemma \ref{lem:difference} implies that for $q\in\mc P_d(G,\ab)$ and $R_hq(g):=q(gh)$, we have
\begin{equation}\label{eq:right-difference}
 R_hq-q\in\mc P_{d-s}(G,\ab).
\end{equation}

\section{Spaces of polynomials on nilspaces}\label{sec:polynomial-spaces}
\noindent The proof of our main theorem relies on the fact that, on polynomial spaces of the form \eqref{eq:polynomials}, there is a nilspace structure that is induced by the degree filtration.
The core properties of this nilspace structure are obtained in this section, by describing the filtration itself in important special cases (see Lemma \ref{lem:integral-polynomials}), and by detailing the associated cube sets (see Lemma \ref{lem:families}).

\subsection{Pure translation lattices as free nilspaces}\hfill\\
The next lemma is the discrete case of \cite[Proposition 3.56]{CGSS-via}, which identifies as free nilspaces those filtered Lie groups having torsion-free, compactly generated successive quotients. We include the splitting argument to make its use here explicit.
\begin{lemma}\label{lem:group-free}
Let $F=\prod_{i=1}^r\mc D_i(\ab_i)$ be a free nilspace and let $\Gamma\le\tran(F)$ be a pure translation lattice. Equip $\Gamma$ with the discrete topology and filtration $\Gamma_i=\Gamma\cap\tran_i(F)$. If $L_i=\eta_i(\Gamma_i)$ has rank $r_i$, then the group nilspace $\Gamma$ is isomorphic to the discrete free nilspace $\prod_{i=1}^r\mc D_i(\mb{Z}^{r_i})$.
\end{lemma}
\begin{proof}
By pureness of $\Gamma$, the homomorphism $\eta_i:\Gamma_i\to L_i$ has kernel equal to $\Gamma_{i+1}$. Thus, we have a short exact sequence
\[
1\longrightarrow\Gamma_{i+1}\longrightarrow\Gamma_i
   \longrightarrow L_i\longrightarrow0
\]
and we have $L_i = \eta_i(\Gamma_i) = \eta_i(\Gamma)\cap\tran_i(F_i)$. Each $L_i$ is free abelian of finite rank by Lemma~\ref{lem:lattices}. These exact sequences also show that $\Gamma$ is finitely generated. Regarding $\Gamma$ as a group nilspace with filtration $\Gamma_\bullet=(\Gamma_i)_{i\ge 0}$, its structure groups are $\Gamma_i/\Gamma_{i+1}\cong L_i$. Starting from the 0-step factor of the group nilspace $\Gamma$, apply \cite[Theorem 4.1]{CGSS-doucos} to each canonical extension and note that an extension of a free nilspace in degree $i$ by the abelian Lie group $L_i$ splits. Hence, induction yields the stated nilspace isomorphism. 
\end{proof}
Note that the isomorphism in Lemma~\ref{lem:group-free} need not be a group homomorphism. Nevertheless, its existence lets us apply lifting results for free nilspaces to $\Gamma$ without choosing group coordinates, as we will see shortly.

\begin{lemma}[Integral polynomial groups]\label{lem:integral-polynomials}
Let $S$ be a discrete free nilspace\footnote{That is, $S$ is a finite product of nilspaces $\mc D_i(\mb{Z}^{r_i})$.}, and let $\Lambda$ be a finitely generated free abelian group. For every integer $d\ge0$, the group $\mc{P}_d(S,\Lambda)$ is finitely generated and free abelian. For every $d\geq 0$ we have $\mc P_{d-1}(S,\Lambda)\le\mc P_d(S,\Lambda)$ and the quotient group $\mc P_d(S,\Lambda)/\mc P_{d-1}(S,\Lambda)$ is free abelian. For $\wt{\ab}=\mb{Z}^a\times\mb{R}^b$, the groups $\mc P_d(S,\wt{\ab})$ are finite-dimensional abelian Lie groups, with the compact-open topology.
\end{lemma}

\begin{proof}
The integral polynomial description in \cite[Lemma 3.5]{CGSS-doucos} expresses a lattice-valued polynomial as a finite sum of products of binomial polynomials, with coefficients in $\Lambda$ and total weighted degree at most $d$. This proves that $\mc{P}_d(S,\Lambda)$ is finitely generated and free. The  subgroup $\mc{P}_{d-1}(S,\Lambda)$ is saturated in the sense that if $nq\in\mc{P}_{d-1}(S,\Lambda)$ for a nonzero integer $n$, then $q\in\mc{P}_{d-1}(S,\Lambda)$, because the alternating-sum equations defining $\mc{P}_{d-1}(S,\Lambda)$ hold for $q$, since multiplication by $n$ is injective on $\Lambda$. The quotient is therefore torsion-free and finitely generated, hence free. The same polynomial description, together with \cite[Lemma 3.6]{CGSS-doucos}, treats $\wt{\ab}$. Coefficients are recovered from finitely many evaluations by finite differences, so the coefficient topology agrees with the compact-open topology.
\end{proof}
The polynomial description is used only to establish these finiteness and topology statements. All function spaces and actions below are defined without a choice of polynomial basis.

\subsection{Cubes of polynomials}\hfill\\
We will describe points by functions on a discrete free nilspace and we will use the following criteria for when a configuration of such functions is a cube. This result is valid more generally than for free nilspaces, and is purely algebraic (not involving any topology).

\begin{lemma}[Polynomial families on a cube]\label{lem:families}
Let $S$ be a nilspace, let $\ab$ be an abelian group, and let $k,n\ge0$. Let $M$ be a subgroup of $\mc{P}_{k-1}(S,\ab)$, equipped with the filtration $M_\bullet$ consisting of the subgroups
\[
 M_0=M_1=M,\qquad M_i=M\cap\mc P_{k-i}(S,\ab)\quad(i\ge1).
\]
For functions $f_v\in M$, indexed by $v\in\cube{n}$, the following conditions are equivalent:\setlength{\leftmargini}{0.7cm}
\begin{enumerate}
\item The configuration $v\mapsto f_v$ is a cube in $\cu^n(M_\bullet)$.
\item For every $m\ge0$ and $d\in\cu^m(S)$, the configuration $(v,w)\mapsto f_v(d(w))$ is a cube in $\cu^{n+m}(\mc{D}_k(\ab))$.
\item The evaluation map $(v,y)\mapsto f_v(y)$ from $\cube{n}\times S$ to $\mc D_k(\ab)$ preserves cubes, where $\cube{n}$ carries the cube structure given by discrete-cube morphisms.
\end{enumerate}
\end{lemma}
Note that the equivalence of (ii) and (iii) also holds for arbitrary functions $f_v:S\to\ab$, without the assumption $f_v\in M$.

\begin{proof}
Let us prove that (ii) $\implies$ (i) $\implies$ (iii)
$\implies$ (ii). For an element $a$ of an abelian group and
a face $F\subset \{0,1\}^n$, write $a^F$ for the corresponding \emph{face-group} element, namely the map on $\{0,1\}^n$ equal to $a$
on $F$ and zero elsewhere.

Suppose that (ii) holds. By Lemma 2.2.5 from \cite{Cand:Notes1}, letting $F_0,\ldots,F_{2^n-1}$ be the upper faces
of $\cube{n}$, ordered as in that lemma,
and letting $s_i=\operatorname{codim}(F_i)$, we have a unique expansion of our configuration as a sum of face-group elements, namely $f_v=\sum_i a_i^{F_i}(v)$, with $a_i\in M$, where each coefficient $a_i$
is an integer linear combination of the functions $f_v$.\footnote{The existence of such an expansion can be proved by evaluating the map $v\mapsto f_v$ in all its vertices as similarly as in the proof of \cite[Lemma 2.2.5]{Cand:Notes1}.}
We must show that $a_i\in M_{s_i}$.
Fix $m\ge0$ and $\q\in\cu^m(S)$. Let $F_j'$ be the upper faces of $\cube{m}$, let $t_j=\operatorname{codim}(F'_j)$, and expand
$a_i\co \q=\sum_j b_{ij}^{F'_j}$, with $b_{ij}\in\ab$.
We then have the upper-face expansion
$f_v(\q(v'))=\sum_{i,j}b_{ij}^{F_i\times F'_j}(v,v')$.
The product faces are precisely the upper faces of
$\cube{n+m}$, and their codimensions are $s_i+t_j$.
By (ii) and \cite[Lemma 2.2.5]{Cand:Notes1}, we have $b_{ij}=0$ whenever $s_i+t_j>k$.
For $s_i\le k$, this means that $a_i\co \q$ is a cube
of $\mc D_{k-s_i}(\ab)$. Since $m,\q$ were arbitrary,
$a_i\in\mc P_{k-s_i}(S,\ab)$.
For $s_i>k$, taking $m=0$ and evaluating at each point of
$S$ gives $a_i=0$. Together with $a_i\in M$, these
conclusions give $a_i\in M_{s_i}$, proving (i).

Now suppose that (i) holds. Consider a cube $(\phi,\q)$ of dimension $r$
in $\cube{n}\times S$, where
$\phi:\cube{r}\to\cube{n}$ is a discrete-cube morphism
and $\q\in\cu^r(S)$. The upper-face expansion above now
has $a_i\in M_{s_i}$. In particular, $a_i=0$ for $s_i>k$,
and $a_i\co \q$ is a cube of
$\mc D_{k-s_i}(\ab)$ for $s_i\le k$.
Expand each $a_i\co \q$ into upper-face elements.
Its nonzero coefficients occur only on faces of
codimension at most $k-s_i$. Therefore the corresponding
product faces in the expansion of
$(v,v')\mapsto f_v(\q(v'))$ have codimension at most $k$.
By \cite[Lemma 2.2.5]{Cand:Notes1}, this configuration
is a cube of $\mc D_k(\ab)$.
Precomposing with the discrete-cube morphism
$v'\mapsto(\phi(v'),v')$ shows that
$v'\mapsto f_{\phi(v')}(\q(v'))$ is also a cube.
This proves (iii).

Finally, assume (iii) and let $\q\in\cu^m(S)$.
The map $(v,v')\mapsto(v,\q(v'))$ is an $(n+m)$-cube
in $\cube{n}\times S$. Applying the evaluation map,
which preserves cubes by (iii), gives
$(v,v')\mapsto f_v(\q(v'))\in
\cu^{n+m}(\mc D_k(\ab))$.
This proves (ii).
\end{proof}

\begin{corollary}\label{cor:families}
Let $S,\ns$ be nilspaces, let $\ab$ be an abelian group, and let $k\ge 0$. Let $P:\ns\times S\to \mc{D}_k(\ab)$ be a morphism. Then, for every $x\in\ns$ the map $g\mapsto P(x,g)$ is in $\mc P_k(S,\ab)$ and for all $x,x'\in\ns$, we have
\begin{equation}\label{eq:slice-difference}
 P(x',\cdot)-P(x,\cdot)\in\mc P_{k-1}(S,\ab).
\end{equation}
\end{corollary}
\begin{proof}
Given a cube $\q\in \cu^n(S)$, if we let $x^{\{0,1\}^n}\in\cu^n(\ns)$ be the constant cube, we have that $P(x,\cdot)\co\q=P^{\{0,1\}^n}(x^{\{0,1\}^n},\q)\in \cu^n(\mc{D}_k(\ab))$ and thus $P(x,\cdot)$ is in $\mc P_k(S,\ab)$. To prove \eqref{eq:slice-difference}, note that since every pair $(x,x')$ is a 1-cube,  for every $\q'\in\cu^m(S)$ the two slices $P(x,\q'(\cdot))$ and $P(x',\q'(\cdot))$ form an $(m+1)$-cube of $\mc D_k(\ab)$. Their difference is a cube of $\mc D_{k-1}(\ab)$ by \cite[Lemma 3.3.37]{Cand:Notes1}.
\end{proof}

\begin{remark}\label{rem:more-general-3.4}
There is a more precise description of the property that the map $P:\ns\times S\to\mc{D}_k(\ab)$ is a morphism. We do not need this in the paper's main arguments, but it is worth mentioning. Fix any $x_0\in\ns$ and define the map $\Psi$ on $\ns$ as $x\mapsto P(x,\cdot)$.  Then, since $P(x,\cdot)=P(x_0,\cdot)+P(x,\cdot)-P(x_0,\cdot)$, we have by \eqref{eq:slice-difference} that $\Psi$ takes values not just in $\mc{P}_k(S,\ab)$ but more precisely in $P(x_0,\cdot)+\mc{P}_{k-1}(S,\ab)$ (note that it follows from \eqref{eq:slice-difference} that this coset is independent of the choice of $x_0$). The more precise description is that then $\Psi$ is also a morphism if we equip $P(x_0,\cdot)+\mc{P}_{k-1}(S,\ab)$ with the standard nilspace structure on $\mc{P}_{k-1}(S,\ab)$ (from Remark \ref{rem:polyfiltra}) transported by the shift by $P(x_0,\cdot)$. This is also more precise in that there is a converse: given a morphism $\Psi:\ns\to \sigma+\mc{P}_{k-1}(S,\ab)$ where $\sigma\in \mc P_k(S,\ab)$, the induced map $\ns\times S\to \mc{D}_k(\ab)$, sending $(x,s)$ to $\Psi(x)(s)$, is a morphism. We omit the details. 
\end{remark}

\section{Pure presentations and toral embeddings}\label{sec:pure-presentations}

\noindent In this section we establish the following result, which is a central ingredient in our proof of Theorem \ref{thm:embedding-intro}.

\begin{theorem}[Pure presentation of \textsc{cfr} nilspaces]\label{thm:pure-presentation}
Every $k$-step \textsc{cfr} nilspace $\ns$ has a presentation $\ns\cong F/\Gamma$, where $F$ is a $k$-step free nilspace and $\Gamma\le\tran(F)$ is a finitely generated pure translation lattice.
\end{theorem}

\subsection{The original extension as a space of point-function pairs}\label{sec:orig-exp-pairs}\hfill\\
We prove Theorem~\ref{thm:pure-presentation} by induction on $k$. For $k=0$ take a point and the trivial group. For $k\ge1$, by induction there is a  presentation
\[
 \ns_{k-1}\cong F/\Gamma
\]
where $F$ is a $(k-1)$-step free nilspace and $\Gamma$ is a pure translation lattice in $\tran(F)$. Let $\varphi:F\to\ns_{k-1}$ denote the quotient fibration. Let $\Gamma_i=\Gamma\cap\tran_i(F)$, so $\Gamma_k=\{1\}$. By \cite[Proposition 5.42]{CGSS-doucos}, the pureness of $\Gamma$ is equivalent to the action of $\Gamma$ on $F$ being free and fiber-transitive. By Lemma~\ref{lem:group-free}, $\Gamma$ is a discrete free nilspace when equipped with its Host--Kra cubes for the filtration $(\Gamma_i)_{i\geq 0}$.

Let $\ab=\ab_k(\ns)$. We now consider the following fiber-product
\[
 \nss=\{(t,x)\in F\times\ns:\varphi(t)=\pi_{k-1}(x)\},
\]
with the topology and cubes inherited from $F\times\ns$. By \cite[Proposition A.20]{CGSS-p-hom}, $\nss$ is a $k$-step nilspace. As $F$ is $(k-1)$-step and $F\to \ns_{k-1}$ is a fibration, the projection $\nss\to F$ makes this a degree-$k$ continuous extension of $F$ by $\ab$ as per \cite[Definition 2.26]{CGSS-doucos}. Note also that both coordinate projections are continuous fibrations. The splitting theorem \cite[Theorem 4.1]{CGSS-doucos} gives a $\ab$-equivariant topological nilspace isomorphism
\[
 \nss\cong F\times\mc D_k(\ab).
\]
On the fiber-product $\nss$, an element $h\in\Gamma$ acts as the map $(t,x)\mapsto (ht,x)$. Note that this map is well-defined (since $\varphi(t)=\varphi(ht)$ by definition). Moreover, if $h\in\Gamma_s$, its action on a codimension-$s$ face preserves the $F$-cube, and leaves the $\ns$-cube unchanged, so $h$ acts on $\nss$ as a translation of height at least $s$. 

Now define the function $\rho:\Gamma\times F\to \ab$ by $h\cdot(t,0)=(ht,\rho(h,t))$ (in particular $\rho(1,t)=0$). Note that $\rho$ is automatically continuous because translations act continuously and $\Gamma$ is a discrete topological space. Lemma~\ref{lem:action} (applied with $G=\Gamma$ and $\nss$ our split extension $F\times\mc D_k(\ab)$ here), followed by projection to the $\ab$ component, shows that $\rho$ is a continuous morphism $\Gamma\times F\to\mc{D}_k(\ab)$. The fact that the action of $\Gamma$ on $\nss$ commutes with the action of $\ab$ implies that, on the chosen split version $F\times \mc{D}_k(\ab)$, the action of $\Gamma$ is of the form\footnote{Indeed, every point $(t,z)$ is obtained from $(t,0)$ by the fiber translation $T_z$. Since $h$ commutes with these translations, we have $h\cdot(t,z)=h\cdot T_z(t,0)=T_z\bigl(h\cdot(t,0)\bigr)=T_z(ht,\rho(h,t))=(ht,\rho(h,t)+z)$.}
\begin{equation}\label{eq:rho-action}
 h\cdot(t,z)=(ht,z+\rho(h,t)).
\end{equation}
This action rule implies the following cocycle identity, by a standard argument:\footnote{We have by \eqref{eq:rho-action} directly $(gh)\cdot (t,z)=((gh)t,z+\rho(gh,t))$, and on the other hand by the action rule combined with \eqref{eq:rho-action} we have $(gh)\cdot (t,z)=g\cdot(h\cdot (t,z))=g\cdot(ht,z+ \rho(h,t))=(ght,z+ \rho(h,t)+\rho(g,ht))$. Equating the second components of these two terms yields the identity.}
\begin{equation}\label{eq:rho-cocycle}
 \rho(gh,t)=\rho(h,t)+\rho(g,ht).
\end{equation}
The projection $\nss\to\ns$ has precisely the $\Gamma$-orbits as its fibers, since two first coordinates above the same $x$ lie in the same fiber of $\varphi$.

The function
\begin{equation}\label{eq:orbfn}
 q_{t,z}:\Gamma \to \ab,\, g\mapsto z+\rho(g,t)
\end{equation}
records the second coordinate of every $\Gamma$-translate of $(t,z)$. We will refer to $q_{t,z} $ as an \emph{orbit function} (the full orbit map would be $g\mapsto g\cdot(t,z)$, which records both coordinates; the orbit function  records only the second coordinate of points in this $\Gamma$-orbit, using the chosen splitting).

Note that by Corollary \ref{cor:families}, the function $g\mapsto \rho(g,t)$ is a polynomial of degree at most $k$. The original second coordinate can be recovered from the orbit function $q_{t,z}$, because $q_{t,z}(1)=z$. Thus we can identify the split extension with the set
\begin{equation}\label{eq:original-pairs}
 \mc E=\{(t,q_{t,z}):t\in F,\ z\in\ab\}
 =\{(t,z+\rho(\,\cdot\,,t)):t\in F,\ z\in\ab\}.
\end{equation}
The bijection $(t,z)\mapsto(t,q_{t,z})$, with inverse $(t,q)\mapsto(t,q(1))$, transports the topology and cubes of $F\times\mc D_k(\ab)$ to $\mc E$. The topology is also the subspace topology from $F\times\ab^\Gamma$, since every evaluation is continuous and the inverse uses evaluation at $1$.

In this space of pairs $\mc{E}$, the cocycle action becomes reindexing of the orbit function. Indeed, defining $R_hq(g)=q(gh)$, we have by \eqref{eq:rho-cocycle} that
\begin{equation}\label{eq:original-reindexing}
 R_hq_{t,z}(g)
 =z+\rho(gh,t)
 =z+\rho(h,t)+\rho(g,ht)
 =q_{ht,z+\rho(h,t)}(g).
\end{equation}
Hence, the original action of $\Gamma$ on the split fiber-product $\nss$ is re-expressed simply as $(t,q)\mapsto(ht,R_hq)$. The identity $R_h\co R_\ell=R_{h\ell}$ shows directly that these transformations of pairs satisfy the multiplication law of $\Gamma$, and therefore every relation in $\Gamma$.

\smallskip

The nilspace $\nss$, when viewed as the split extension $F\times \mc{D}_k(\ab)$, is useful on the way to our goal of obtaining a free nilspace $\mc{F}$ (defined further below in \eqref{eq:profile-space}) for our desired pure presentation of $\ns$. However, the main remaining issues include that $F\times \mc{D}_k(\ab)$ is clearly not a free nilspace. Presenting it as the space $\mc{E}$ is a key step towards defining $\mc{F}$. The advantage of $\mc{E}$ lies in its expression as a space of pairs involving the cocycle $\rho$ in the second component. Indeed, this expression affords the possibility to define our final construction in a natural way as a similar space of pairs, but where the second component involves lifts of the orbit functions  to polynomial functions of appropriate degree taking values in an abelian Lie group
$\wt{\ab}=\mb{Z}^a\times\mb{R}^b$ covering $\ab$,
so that $\mc{D}_k(\wt{\ab})$ is a free nilspace. Formalizing this construction is the purpose of the next subsection.

\subsection{Polynomial lifts and normalized lattice-valued functions}\label{sec:with-example}\hfill\\
As $\ab$ is a compact abelian Lie group, we may choose a covering homomorphism $p$ for it, associated with a decomposition of $\ab$ into a product of finite cyclic groups and a torus. More precisely, supposing that $\ab=B\times \mb{T}^b$ where $B$ is a finite abelian group of rank $a$, we have a short exact sequence
\begin{equation}\label{eq:top-cover}
0\longrightarrow\Lambda\longrightarrow\wt{\ab}
 \overset p\longrightarrow\ab\longrightarrow0,
 \qquad \wt{\ab}=\mb{Z}^a\times\mb{R}^b,
\end{equation}
where the kernel $\Lambda$ of $p$ is discrete and cocompact, and hence free abelian of finite rank by Lemma~\ref{lem:lattices}. Since $\Gamma\times F$ is free, by \cite[Theorem 3.11]{CGSS-doucos}, the map $\rho$ can be lifted to a continuous morphism $P:\Gamma\times F\to\mc{D}_k(\wt{\ab})$. Subtracting $P(1,t)$ from $P(g,t)$ preserves the morphism property and its image under $p$, because $P(1,t)\in\Lambda$. We may therefore require
\begin{equation}\label{eq:joint-lift}
 p\co P=\rho,\qquad P(1,t)=0.
\end{equation}
Let us introduce the shorthand notation $P_t(g):=P(g,t)$. Equation \eqref{eq:slice-difference} then gives
\begin{equation}\label{eq:lift-difference}
 P_{t'}-P_t\in\mc P_{k-1}(\Gamma,\wt{\ab}) \text{ for all }t,t'\in F.
\end{equation}
The idea now is to try to see the action of $\Gamma$ on pairs $(t,P_t+z)$ where $z\in\wt{\ab}$ as before, i.e. $(t,P_t+z)\mapsto (ht,R_hP_t+z)$. However, the problem is that this is not well-defined in general, as the next example shows.

\begin{example}\label{ex:problematic-3-step}
Let $F=\mc{D}_1(\mb{Z}^2)$, let $\wt{\ab}=\mb{Z}$, let $\Gamma:=\langle h_1,h_2\rangle$ where $h_1(x,y)=(x+2,y)$ and $h_2(x,y)=(x,y+2)$, and let $P:\Gamma\times F \to \mc{D}_3(\mb{Z})$ be defined by $P(h_1^nh_2^m,(x,y)):=mx^2$ (in particular $P_{(x,y)}(n,m)=mx^2$). It can be checked that, letting $\Lambda:=4\mb{Z}$, this is a valid lift of a $\rho$ which defines a degree-3 extension by $\mb{Z}_4$ of $F/\Gamma\cong\mc{D}_1(\mb{Z}^2_2)$. If we apply the previous shift formula $(t,P_t+z)\mapsto (ht,R_hP_t+z)$ to $((x,y),P_{(x,y)}+z)$ with $h_1$, we get $((x+2,y),P_{(x,y)}+z)$. However, this is not equal to $((x+2,y),P_{(x+2,y)}+z')$ for any constant $z'\in \mb{Z}$, as the difference is the non-constant morphism $(n,m)\mapsto m(4x+4)$.
\end{example}
We resolve this issue by allowing a larger family of functions appearing as the second component of the pairs $(t,P_t+z)$. This will lead to the introduction of the larger space $\mc{F}'$ defined in the next subsection. To close the present subsection we introduce the following objects, which will be useful in what follows.

Let
\begin{equation}\label{eq:auxiliary-filtration}
 V=\mc P_{k-1}(\Gamma,\Lambda),\qquad
 V_i=\mc P_{k-i}(\Gamma,\Lambda)\ (i\ge1),\qquad V_0=V_1,
\end{equation}
and define the subgroup of normalized polynomials
\begin{equation}\label{eq:normalized-U}
 U=\{u\in V:u(1)=0\},
\end{equation}
with filtration $(U_i=U\cap V_i)_{i\ge0}$. Note that $U_0=U_1$ and $U_k=0$ (thus the filtration on $U$ is of degree $k-1$).

Note that the subgroup of constants of $V$ is equal to $V_k=\Lambda$, and that we may decompose any $w\in V$ uniquely as
\[
 w=w(1)+(w-w(1))\; \text{ where }\;w(1)\in \Lambda \;\text{ and }\;w-w(1)\in U.
\]
\noindent We will refer to the map $V\to U$, $w\mapsto w-w(1)$ as the \emph{normalization}. 
For $1\le i<k$, note that $v(1)\in V_{i+1}$ and that $v-v(1)\in U_i$ and we have the isomorphism
\begin{equation}\label{eq:isou}
U_i/U_{i+1}\cong V_i/V_{i+1}.
\end{equation}
Indeed, the homomorphism $U_i\to V_i/V_{i+1}$, 
$u\mapsto u+V_{i+1}$ has kernel
$U_i\cap V_{i+1}=U_{i+1}$, and is surjective because for $v\in V_i$, the constant map with value $v(1)$ lies in $V_{i+1}$ (as $i<k$), while $v-v(1)\in U_i$, so that every class
modulo $V_{i+1}$ has a normalized representative.

The following result will be used to ensure that the main construction in the next subsection (the nilspace $\mc{F}$) is a free nilspace.
\begin{lemma}\label{lem:U-is-free}
The group $U$, equipped with the filtration $(U_i)_{i\ge 0}$ and the discrete topology, is a free nilspace of step at most $k-1$.    
\end{lemma}
\begin{proof}
The quotients $U_i/U_{i+1}$ are finitely generated free abelian groups (e.g.\ by \eqref{eq:isou} and Lemma~\ref{lem:integral-polynomials}). We can therefore choose a complement $W_i$ of $U_{i+1}$ in $U_i$, successively for $i=k-1,\ldots,1$. The resulting filtered-group decomposition identifies $U$, with its discrete topology, with the free nilspace $\prod_{i=1}^{k-1}\mc D_i(W_i)$ by \cite[Theorem 4.1]{CGSS-doucos}.
\end{proof}

\subsection{The ambient action and the invariant set of point-function pairs}\label{subsec:invaset}\hfill\\
Consider first the set of all point-function pairs
\begin{equation}\label{eq:ambient-pairs}
 \mc F'=\{(t,f):t\in F,\ f:\Gamma\to\wt{\ab}\}=F\times \wt{\ab}^\Gamma.
\end{equation}
We equip this set with the product topology from $F\times\wt{\ab}^\Gamma$ (which is second-countable, since $\Gamma$ is countable by Lemma \ref{lem:group-free}). We do not need to define a cube structure on $\mc{F}'$. Instead, we first define a semidirect product group $V\rtimes\Gamma$, which acts on $\mc{F}'$, and which will play a key role eventually as our new translation lattice. This requires
specifying an action of $\Gamma$ by automorphisms of the
additive group $V$. 

For $h\in\Gamma$, recall that the reindexing operator $R_h$ on $\wt{\ab}^\Gamma$ is defined by $(R_hf)(g)=f(gh)$.
\begin{lemma}
The set $V\times\Gamma$ with multiplication
\begin{equation}\label{eq:profile-semidirect-product}
 (w,h)(w',h')=(w+R_hw',hh'),
\end{equation}
is a semidirect product group $V\rtimes\Gamma$. Moreover, this group acts on $\mc{F}'$ by homeomorphisms, as follows: for $(w,h)\in V\rtimes\Gamma$ and
$(t,f)\in\mc F'$, we let
\begin{equation}\label{eq:profile-group-action}
 (w,h)\cdot(t,f):=(ht,R_hf+w).
\end{equation}
\end{lemma}
\begin{proof}
We first prove that \eqref{eq:profile-semidirect-product} indeed defines a binary operation on $V\times \Gamma$, i.e.\ that the right side is in $V\times \Gamma$. For this it suffices to check that $R_h$ maps $V$ to itself. Right multiplication $g\mapsto gh$ preserves the Host--Kra cubes
of $\Gamma$, indeed if $\q$ is such a cube, then $v\mapsto \q(v)h$ is the pointwise product of $\q$ and the
constant cube with value $h$, and is therefore a cube.
Thus, for $w\in V=\mc P_{k-1}(\Gamma,\Lambda)$,
the shifted map $R_hw$ is again a $\Lambda$-valued
polynomial of degree at most $k-1$. Hence $R_h(V)\subseteq V$. (Alternatively, check that $g\mapsto gh$ is a translation and apply Lemma \ref{lem:difference} to see that $R_h(V)\subseteq V$.)
Moreover, $R_h$ is additive and has inverse $R_{h^{-1}}$,
so it is an automorphism of $V$.

To obtain an action of $\Gamma$ on $V$, we also need these
automorphisms to respect multiplication in $\Gamma$.
For $h,\ell,g\in\Gamma$ and $w\in V$, we have
$(R_hR_\ell w)(g)=w((gh)\ell)=w(g(h\ell))$.
Therefore $R_hR_\ell=R_{h\ell}$, and $R_1$ is the identity.

The automorphism and composition properties just proved
ensure that this multiplication is associative.
Its identity is $(0,1)$, and its inverses are
$(w,h)^{-1}=(-R_{h^{-1}}w,h^{-1})$.

We next verify the action of this group $V\rtimes\Gamma$ on the ambient space $\mc{F}'$. 
The second coordinate is the function
$g\mapsto f(gh)+w(g)$, where we regard
$\Lambda$ as a subgroup of $\wt{\ab}$.

We verify that this formula respects the group
multiplication, as required for a left action.
Using additivity of $R_h$ and the identity
$R_hR_{h'}=R_{hh'}$, we obtain
\begin{align*}
(w,h)\cdot\bigl((w',h')\cdot(t,f)\bigr)&=\bigl(h(h't),R_h(R_{h'}f+w')+w\bigr)\\
&=\bigl((hh')t,R_{hh'}f+R_hw'+w\bigr) =\bigl((w,h)(w',h')\bigr)\cdot(t,f).
\end{align*}
The identity $(0,1)$ fixes every pair, so
\eqref{eq:profile-group-action} defines a left action.
Each transformation is a homeomorphism: the map
$t\mapsto ht$ is a homeomorphism of $F$, reindexing
permutes the factors of $\wt{\ab}^{\Gamma}$, and addition
by the fixed function $w$ is a homeomorphism.
The inverse transformation is the action of
$(w,h)^{-1}$.
\end{proof}

\begin{remark}[Wreath product]
A related construction is the unrestricted regular \emph{wreath product} of $\wt{\ab}$ by $\Gamma$. This is defined as $W=\wt{\ab}^{\Gamma}\rtimes_R\Gamma$, where $\wt{\ab}^{\Gamma}$ consists of all functions $\Gamma\to\wt{\ab}$, added pointwise, and $R_hf(g)=f(gh)$. Note that this group acts on $\mc{F}'$ by $(f',h)\cdot(t,f)=(ht,R_hf+f')$. Since $R_h(V)=V$, the group $V\rtimes\Gamma$ is a subgroup of $W$.
\end{remark}

While the space $\mc F'$ is large enough to accommodate the action of $V\rtimes\Gamma$, it is too large to represent the \textsc{cfr} nilspace $\ns$ adequately. In particular, it lacks a suitable free nilspace structure. 
To address this issue, we will restrict the action of $V\rtimes\Gamma$ to a subset $\mc{F}\subset \mc{F}'$, defined below. The space $\mc F$ will be equipped with a $k$-step free nilspace structure. Our main task will be to prove that $V\rtimes\Gamma$ acts on it freely and fiber-transitively, and finally that $\mc F/(V\rtimes \Gamma)\cong \ns$.

Let us first establish the freeness of the action, which actually holds on all of $\mc{F}'$.
\begin{lemma}
The action of $V\rtimes\Gamma$ on $\mc F'$ is free.
\end{lemma}
\begin{proof}
If $(w,h)\cdot(t,f)=(t,f)$, equality of the first
coordinates gives $ht=t$. Freeness of the lower
$\Gamma$-action on $F$ implies $h=1$.
Equality of the second coordinates then gives $f+w=f$,
hence $w=0$. Therefore only the identity group element
can fix a pair. 
\end{proof}
Any invariant subset of $\mc F'$ will inherit this free action.

We now define the subset $\mc F$ that will be our new free nilspace:
\begin{equation}\label{eq:profile-space}
 \mc F=\{(t,\wt{z}+P_t+u):t\in F,\ \wt{z}\in\wt{\ab},\ u\in U\}
 \subseteq\mc F'.
\end{equation}
\noindent Note that reduction modulo $\Lambda$ sends a function $f(g)=\wt{z}+P_t(g)+u(g)$ in the second coordinate to
\begin{equation}\label{eq:profile-reduction}
p(f(g))=p(\wt{z})+\rho(g,t)=q_{t,p(\wt{z})}(g).
\end{equation}
Thus, these point-function pairs $(t,\wt{z}+P_t+u)$ lift precisely the original pairs \eqref{eq:original-pairs}, with the normalized corrections $u$ recording the auxiliary choices. Conversely, every original pair in the split fiber-bundle $\nss$ (or $\mc{E}$) has a lift, since $p$ is surjective.

\begin{lemma}\label{lem:profile-invariant}
The set $\mc F$ is invariant under the action of $V\rtimes\Gamma$ on $\mc F'$.
\end{lemma}
For the upcoming proof, it is useful to separate
this action into reindexing and addition. Denote the
actions of $(0,h)$ and $(w,1)$ by $A_h$ and $T_w$,
respectively:
\begin{equation}\label{eq:profile-actions}
 A_h(t,f)=(ht,R_hf),\qquad
 T_w(t,f)=(t,f+w).
\end{equation}
Since $(w,h)=(w,1)(0,h)$, every  element of $V\rtimes\Gamma$ acts as
$T_wA_h$. Thus invariance under all $A_h$ and $T_w$
will imply invariance under the whole group.
For computations with these transformations, the
semidirect-product multiplication gives
\begin{equation}\label{eq:profile-relations}
 A_hA_\ell=A_{h\ell},\qquad
 A_hT_wA_h^{-1}=T_{R_hw}.
\end{equation}

\begin{proof}[Proof of Lemma \ref{lem:profile-invariant}]
Fix any point $(t,f)\in\mc{F}$, where $f=\wt{z}+P_t+u$. 

We first show that for any $w\in V$, we have $T_w(t,f)\in \mc{F}$. Note that $f+w=(\wt{z}+w(1))+P_t+(u+w-w(1))$, and the last summand belongs to $U$, so indeed $T_w(t,f)\in \mc{F}$.

Next we prove that for every $h\in\Gamma$, the pair $A_h(t,f)=(ht,R_hf)$ is in $\mc{F}$. For this, we must express its function coordinate in the form $R_hf=\wt{z}'+P_{ht}+u'$ with $u'\in U$. The only possible choices for $\wt{z}'$ and $u'$ are 
\begin{equation}\label{eq:reindexed-coordinates}
 \wt{z}'=f(h)\quad\text{and}\quad u'=R_hf-f(h)-P_{ht}.
\end{equation}
Indeed, if a desired decomposition $R_hf=\wt{z}'+P_{ht}+u'$ exists, then evaluating it at the identity $1\in\Gamma$ gives $f(h)=(R_hf)(1)=\wt{z}'+P_{ht}(1)+u'(1)=\wt{z}'$, because $P_{ht}(1)=0$ and $u'(1)=0$. Therefore $\wt{z}'$ must equal $f(h)$. Subtracting these terms in the initial decomposition then forces $u'=R_hf-f(h)-P_{ht}$. 

We now verify that $u'\in U$. First note that from \eqref{eq:reindexed-coordinates} we get $u'(1)=0$, because $P_{ht}(1)=0$ and $f(h)=(R_hf)(1)$. Secondly, applying the homomorphism $p$ and using \eqref{eq:profile-reduction}, we have $p(u'(g))=\rho(gh,t)-\rho(h,t)-\rho(g,ht)$ and this is $0$ by the cocycle identity \eqref{eq:rho-cocycle}, so $u'$ is $\Lambda$-valued. Finally, write
\[
 u'=(R_hf-f)+(f-P_{ht})-f(h).
\]
Since $f\in\mc{P}_k(\Gamma,\wt{Z})$, the first term $R_hf-f$ here is in $\mc{P}_{k-1}(\Gamma,\wt{Z})$ by \eqref{eq:right-difference}, with $h\in\Gamma_1$. The second term $f-P_{ht}=\wt{z}+(P_t-P_{ht})+u$ is in $\mc{P}_{k-1}(\Gamma,\wt{Z})$ by \eqref{eq:lift-difference} and since $u\in\mc{P}_{k-1}(\Gamma,\Lambda)$ and $\wt{z}$ is constant. The last term $f(h)$ is constant. Therefore $u'\in\mc P_{k-1}(\Gamma,\Lambda)$ and $u'(1)=0$, as required. We have thus proved that $A_h(t,f)=(ht,\wt{z}'+P_{ht}+u')\in\mc F$. The same argument applies to $-w$ and $h^{-1}$, so the inverses preserve $\mc F$ as well.
\end{proof}
We thus have a free action of $V\rtimes \Gamma$ on $\mc{F}$. This action gives a faithful group of homeomorphisms of $\mc{F}$ with its subspace topology, which will  be our new translation lattice:
\begin{equation}\label{def:Gamma'}
 \Gamma'=\{T_wA_h:w\in V,\ h\in\Gamma\}\cong V\rtimes\Gamma.
\end{equation}
Note that $\Gamma'$ is finitely generated because $V$ and $\Gamma$ are.

\begin{remark}
As noted above, the product topology on $\mc{F}'$ is second-countable. It is not necessarily locally compact (the space can involve countable powers of $\mb{R}$). However, when we restrict this topology to the subset $\mc{F}$, the resulting relative topology \emph{is} locally compact. This can be seen for instance using the homeomorphism $\chi$ presented in the next subsection.
\end{remark}

Our next goals are to equip $\mc F$ with a useful free nilspace structure (in Subsection \ref{subsec:freestruct}), and then to complete the main task of proving that $\Gamma'$ is a pure translation lattice on this free nilspace (Subsection \ref{subsec:mainpfs}).

\subsection{The free nilspace structure on $\mc{F}$ and its cubes}\label{subsec:freestruct}\hfill\\
As seen in the proof of Lemma \ref{lem:profile-invariant}, for each $(t,f)\in\mc{F}$ the parameters $\wt{z}$ and $u$ giving the decomposition $f=\wt{z}+P_t+u$ in \eqref{eq:profile-space} are unique, namely
\begin{equation}\label{eq:profile-inverse}
 \wt{z}=f(1)\quad \text{and}\quad u=f-f(1)-P_t.
\end{equation}
Thus the map
\begin{equation}\label{eq:profile-chart}
 \chi:F\times\mc D_k(\wt{\ab})\times U\longrightarrow\mc F,
 \qquad (t,\wt{z},u)\longmapsto(t,\wt{z}+P_t+u)
\end{equation}
is a bijection. We can then equip $\mc{F}$ with the cube structure transported by $\chi$. Since $U$ is free of step at most $k-1$, this makes $\mc F$ a free nilspace of step at most $k$.

The product topology in \eqref{eq:profile-chart} agrees with the subspace topology already given to $\mc F$. Each evaluation of $\wt{z}+P_t+u$ depends continuously on $(t,\wt{z},u)$, so $\chi$ is continuous. For the inverse, $\wt{z}=f(1)$ is continuous, and each evaluation of $f-f(1)-P_t$ is continuous. By Lemma~\ref{lem:integral-polynomials}, finitely many evaluations recover a bounded-degree polynomial's coefficients. Because these coefficients lie in the discrete lattice $\Lambda$, the pointwise topology on $U$ is discrete. Thus the inverse in \eqref{eq:profile-inverse} is continuous. In particular, the restricted ambient action consists of homeomorphisms for the free nilspace topology.

The following criterion expresses the cubes directly in terms of the pairs. 
\begin{lemma}[Cubes on $\mc{F}$]\label{lem:profile-cubes}
A map $(\q_t,\q_f):\{0,1\}^n\to \mc{F}$ sending $v$ to $(t_v,f_v)$ is an $n$-cube on $\mc{F}$ if and only if $\q_t\in\cu^n(F)$ and the evaluation map $(v,g)\mapsto f_v(g)$ 
from $\cube{n}\times\Gamma$ to $\mc D_k(\wt{\ab})$ preserves cubes. Equivalently, in addition to the first-coordinate condition, for every $m\ge0$ and $\q\in\cu^m(\Gamma)$ the configuration $(v,v')\mapsto f_v(\q(v'))$ is an $(n+m)$-cube.
\end{lemma}
\begin{proof}
Write $f_v=\wt{z}_v+P_{t_v}+u_v$ using \eqref{eq:profile-inverse}. Suppose first that the configuration is a cube. Then $\q_{\wt{z}}$ defined as $\q_{\wt{z}}(v):=\wt{z}_v$ is a cube of $\mc D_k(\wt{\ab})$, $\q_t$ is a cube of $F$, and $\q_u$ defined as $\q_u(v)=u_v$ is a cube of $U$. The evaluation maps of the first two summands preserve cubes, using the joint morphism $P$ for the second. For the third, Lemma~\ref{lem:families}, applied to $M=U\le\mc P_{k-1}(\Gamma,\Lambda)$ with its filtration \eqref{eq:normalized-U}, shows that $(v,g)\mapsto u_v(g)$ preserves cubes in $\mc D_k(\Lambda)$, and hence also in $\mc D_k(\wt{\ab})$. Adding the three evaluation maps proves the required condition.

Conversely, suppose $\q_t$ is a cube and the joint evaluation preserves cubes. Evaluation at $g=1$ gives a cube $\q_{\wt{z}}$ of $\mc D_k(\wt{\ab})$. Subtracting $\wt{z}_v+P_{t_v}(g)$ from the evaluation leaves $u_v(g)$, which also preserves cubes. This evaluation is $\Lambda$-valued. For every $r\ge0$, a configuration
$\q':\cube{r}\to\Lambda$ belongs to $\cu^r(\mc D_k(\Lambda))$
if and only if, viewed as a $\wt{\ab}$-valued configuration, it belongs
to $\cu^r(\mc D_k(\wt{\ab}))$, since the inclusion
$\Lambda\hookrightarrow\wt{\ab}$ preserves the
alternating-sum equations defining these cubes. Lemma~\ref{lem:families}, with the same choice $M=U$, therefore shows that $\q_u$ is a cube of $U$. All three coordinates in \eqref{eq:profile-chart} are therefore cubes. The final equivalent statement follows from the fact that we can fix an arbitrary $t_0\in F$ and regard $f_v=P_{t_0}+f_v'$ with $f_v'\in \mc P_{k-1}(\Gamma,\wt \ab)$ and $P_{t_0}\in\mc P_k(\Gamma,\wt \ab)$ by Corollary \ref{cor:families}. Then, the evaluation map $(v,g)\mapsto f'_v(g)$ preserves cubes if and only if $(v,g)\mapsto f_v(g)$ preserves cubes. Thus, the result follows from the equivalence of (ii) and (iii) in Lemma~\ref{lem:families} applied to $f_v'$.
\end{proof}

\begin{remark}
There are other enlightening viewpoints on the nilspace structure we have chosen on $\mc{F}$ via \eqref{eq:profile-chart}. 

For instance, it is natural to wonder if there is a nilspace structure on the larger space $\mc{F}'$ itself which restricts precisely to the chosen nilspace structure on $\mc{F}$. A simple candidate is the product-nilspace structure  $F\times \prod_{\gamma\in \Gamma} \mc{D}_k(\wt{\ab})$, but the restriction of this structure to $\mc{F}$ gives strictly more cubes than in the chosen structure when $U\neq \{0\}$. However, letting $\mc{H}=\wt{\ab}^\Gamma$ be regarded as an abelian group under pointwise addition, and defining the degree-$k$ filtration $\mc{H}_0=\mc{H}_1=\mc{H}$ and $\mc{H}_i=\mc{P}_{k-i}(\Gamma,\wt{\ab})
$ for $i\ge2$ (using also the convention $\mc{P}_d={0}$ for $d<0$), we can equip $\mc{H}$ with the associated group nilspace structure, and then equip $\mc{F'}=F\times\mc{H}$ with the product nilspace structure. It can be checked that the restriction of this structure to $\mc{F}$ yields exactly the chosen nilspace structure.

Another viewpoint consists in defining $\mc F$ as a fiber-product
of nilspaces, which is a way to endow it directly with the
fiber-product nilspace structure. 
Indeed, fix any $t_0\in F$ and consider the maps
$p^*:(P_{t_0}+\mc P_{k-1}(\Gamma,\wt{\ab}))
\to(\rho_{t_0}+\mc P_{k-1}(\Gamma,\ab))$
given by $f\mapsto p\co f$, and
$\beta:F\times\mc D_k(\ab)
\to(\rho_{t_0}+\mc P_{k-1}(\Gamma,\ab))$
given by $(t,z)\mapsto q_{t,z}(\cdot)=z+\rho(\cdot,t)$.
These two maps are both morphisms\footnote{On each coset $\sigma+M$ below,
with $M$ equipped with the filtration of
Lemma~\ref{lem:families}, a configuration $v\mapsto f_v$ is
a cube precisely when $v\mapsto f_v-\sigma$ is a cube of $M$. This is the intended nilspace structure on $P_{t_0}+\mc P_{k-1}(\Gamma,\wt{\ab})$ relative to which $p^*$ is a morphism.} (for $p^*$ this is clear
since $p$ is a homomorphism, and for $\beta$ we apply
Remark~\ref{rem:more-general-3.4}). As $p^*$ is a fibration
(by polynomial lifting on $\Gamma$ and the upper-face
description of cubes; see
\cite[Theorem 3.11]{CGSS-doucos}), the fiber product
$(F\times\mc D_k(\ab))
\times_{\rho_{t_0}+\mc P_{k-1}(\Gamma,\ab)}
(P_{t_0}+\mc P_{k-1}(\Gamma,\wt{\ab}))$
is a well-defined nilspace (this follows from the proof of
\cite[Proposition A.7]{CGSSspectral}). It can be checked that
this construction is independent of the chosen element
$t_0\in F$, and that the map $((t,z),f)\mapsto(t,f)$ is a nilspace isomorphism onto $\mc{F}$.
\end{remark}

\subsection{Translation heights, fiber-transitivity and pureness}\label{subsec:mainpfs}\hfill\\
\begin{lemma}\label{lem:profile-heights}
For every $s\in [k]$, we have
\[
 h\in\Gamma_s\Rightarrow A_h\in\tran_s(\mc F),\quad \text{and }\quad
 w\in V_s\Rightarrow T_w\in\tran_s(\mc F).
\]
\end{lemma}
\begin{proof}
For a nilspace $W$, we use the notation (used in the proof of Lemma \ref{lem:families}) where given a function $\alpha:W\to W$, and a face $D\subset \{0,1\}^n$, the map $\alpha^D:W^{\{0,1\}^n}\to W^{\{0,1\}^n}$ is defined by $\alpha^D(f)(v)=\alpha(f(v))$ if $v\in D$ and $\alpha^D(f)(v)=f(v)$ otherwise.

Let $\q_{\mf{f}}=(\q_{t},\q_{f})\in \cu^n(\mc F)$, where $n\ge s$, and let $D\subseteq\cube{n}$ be a face of codimension $s$. We have already proved invariance of $\mc F$ and continuity of the actions and their inverses, so it remains to check the cubes after applying each action on $D$.

For $h\in\Gamma_s$, we have that $A_h^D(\q_{\mf{f}})=(h^D\q_{t}, R_{h^D}\q_{f})$. As $h\in\Gamma_s$, by definition of translations on $F$, we have that $h^D\q_{t}\in \cu^n(F)$. To deal with the term $R_{h^D}\q_{f}$, we use the characterization given by Lemma~\ref{lem:profile-cubes}. Let us denote $\q_f(v)$ by $f_v$ for each $v$. Hence, it suffices to prove that, given $\q\in \cu^m(\Gamma)$, the map $(v,v')\mapsto R_{h^D}f_v(\q(v'))$ is in $\cu^{n+m}(\mc{D}_k(\wt{\ab}))$. Now note that $R_{h^D}f_v(\q(v')) = f_v(\q(v')h^D(v))$. Moreover, note that the map $(v,v')\mapsto \q(v')h^D(v)$ is a Host--Kra cube in $\cu^{n+m}(\Gamma_\bullet)$, because it is a product of two elements of this group, namely (the obvious extensions to $\{0,1\}^{n+m}$ of) $\q$ and $h^D$. 

To conclude, apply Lemma \ref{lem:profile-cubes} to the
original cube $\q_{\mf{f}}=(\q_t,\q_f)$.
Recalling that $f_v=\q_f(v)$, this lemma tells us that
the evaluation map $E:\cube{n}\times\Gamma\to
\mc{D}_k(\wt{\ab})$, defined by $E(v,g)=f_v(g)$,
preserves cubes. We apply this morphism $E$ to the map
$(v,v')\mapsto(v,\q(v')h^D(v))$.
This map is an $(n+m)$-cube of $\cube{n}\times\Gamma$. Indeed, its first component is the coordinate projection, and its
second component is the Host--Kra cube just established above.
Its image under $E$ is therefore the $(n+m)$-cube
$(v,v')\mapsto f_v(\q(v')h^D(v))$ that we
needed. Together with the first-coordinate condition
already proved, Lemma~\ref{lem:profile-cubes} shows that
applying $A_h$ on $D$ preserves the original cube.

Let us now prove that if $w\in V_s$ then $T_w\in \tran_s(\mc F)$. We have that $T_w(\mf{f})=T_w(t,f)=(t,f+w)$. We can write $w=w(1)+(w-w(1))$. To see that this is a translation of height $s$, we use the cube structure that we have given to $\mc F$ via $\chi$ from \eqref{eq:profile-chart}. Note that in the coordinates $(t,\wt z,u)\in F\times\mc{D}_k(\wt \ab)\times U$, we have that $T_w(t,\wt z,u)=(t,\wt z+w(1),u+w-w(1))$. As $w(1)\in \Lambda$ is a constant, we can regard it as a height-$k$ translation, which is in $\tran_k(\mc{D}_k(\wt \ab))$. On the other hand, the element $w-w(1)$ is by definition in $U_s$, and as the nilspace $U$ is equipped with the group nilspace structure, adding an element of degree $s$ is always a translation in $\tran_s(U)$. The result follows.
\end{proof}

Taking $s=1$ proves that $\Gamma'\le\tran(\mc F)$. The freeness of the action of $\Gamma'$ has already been established on $\mc F'$. We now check fiber-transitivity.

\begin{lemma}\label{lem:profile-fiber-transitive}
The action of $\Gamma'$ on $\mc{F}$ is fiber-transitive.
\end{lemma}

\begin{proof}
We will prove that, for every $s\in [1,k]$ and $\mf{f}=(t,f)\in \mc{F}$, \begin{equation}\label{eq:transcharac}
\mf{f}'=T_wA_h\mf{f}\text{ and }\pi_{s-1}(\mf{f})=\pi_{s-1}(\mf{f}')\; \Longrightarrow \; h\in\Gamma_s\text{ and }w\in V_s.
\end{equation}
By Lemma~\ref{lem:profile-heights}, these two memberships
will imply that the given transformation $T_wA_h$ \emph{itself}\footnote{The fact that the translation $T_wA_h$  itself has this required height is apparently stronger than fiber-transitivity, but given that the action is already known to be free, this fact is actually equivalent to fiber-transitivity.} has
height at least $s$, which will confirm fiber-transitivity.

We first prove that $h\in\Gamma_s$. Projecting the assumed
equality $\pi_{s-1}(\mf{f})=\pi_{s-1}(\mf{f}')$ to $F$ gives
$\pi_{s-1}(t)=\pi_{s-1}(ht)$. If $s<k$, the available fiber-transitivity of $\Gamma$ on $F$ provides $h_s\in\Gamma_s$ with
$h_st=ht$. Since the action of $\Gamma$ on $F$ is free, we have $h_s=h$.
If $s=k$, then $F$ has step at most $k-1$, so its
$(s-1)$-step factor is $F$ itself. Thus $ht=t$, and
freeness of $\Gamma$ gives $h=1\in\Gamma_k$. We have therefore
proved that $h\in\Gamma_s$ in every case.

Next, we show that $w\in V_s$. Let $\mf{f}''=A_h\mf{f}$.
Since $h\in\Gamma_s$, Lemma~\ref{lem:profile-heights}
gives that $A_h$ has height at least $s$ and so $\pi_{s-1}(\mf{f}'')=\pi_{s-1}(\mf{f})=\pi_{s-1}(\mf{f}')$.
Moreover $\mf{f}''=(ht,R_hf)$ and $\mf{f}'=(ht,R_hf+w)$.

To obtain information about $w$, compare the coordinate $u\in U$ of these two points $\mf{f}''$ and $\mf{f}'$ as per formula \eqref{eq:profile-inverse}. Their first
coordinate is the same element $ht$, so their term $P_{ht}$ in
\eqref{eq:profile-inverse} is the same.
Since $(R_hf+w)(1)=(R_hf)(1)+w(1)$, that formula shows
that their $U$-coordinates differ by $w-w(1)$. On the other hand, the product nilspace isomorphism
\eqref{eq:profile-chart} implies that the $(s-1)$-step factor
of $\mc{F}$ is
$F_{s-1}\times(U/U_s)$. Indeed, the factor $\mc{D}_k(\wt{\ab})$ contributes only
a point because $s-1<k$.
More explicitly, the point $\chi(t,\wt z,u)$ has image
$(\pi_{s-1}(t),u+U_s)$ in this product.
Therefore the equality
$\pi_{s-1}(\mf f'')=\pi_{s-1}(\mf f')$ implies that their
$U$-coordinates have the same image in $U/U_s$.
Since those coordinates differ by $w-w(1)$,
we obtain $w-w(1)\in U_s\subseteq V_s$. The constant function $w(1)$ also belongs to $V_s$, since its value lies in $\Lambda$ and its degree is $0\le k-s$. Hence $w=w(1)+(w-w(1))\in V_s$, which completes the proof.
\end{proof}
We can now complete the proof of the following main result of this subsection.

\begin{proposition}\label{prop:profile-pure}
The group $\Gamma'$ is a pure translation lattice.
\end{proposition}
\begin{proof}
We have already proved that $\Gamma'$ acts freely on $\mc{F}$ and is fiber-transitive, so the pureness of $\Gamma'$ follows by \cite[Proposition 5.42]{CGSS-doucos}.

It remains to prove that $\Gamma'$ is fiber-discrete and fiber-cocompact. For this, it will help to prove that
\begin{equation}\label{eq:pure-filtration}
 \{T_wA_h:w\in V_s,\ h\in\Gamma_s\}
=\Gamma'\cap\tran_s(\mc F)
 \quad\text{for}\quad 1\le s\le k.
\end{equation}
To see this, note that the left-hand set is included in the right-hand set by Lemma \ref{lem:profile-heights}.  For the opposite inclusion, note that if a translation $T_wA_h$ is in $\tran_s(\mc{F})$, then it fixes the entire $(s-1)$-step factor (i.e. it is in $\ker(\eta_{s-1})$), and then applying \eqref{eq:transcharac} at any point $\mf{f}\in \mc{F}$ gives that $h\in \Gamma_s$ and $w\in V_s$, as required.

We now compute the subgroups whose discreteness and
cocompactness are required. Let
$\Gamma'_j=\Gamma'\cap\tran_j(\mc F)$ and let
$K_j=\eta_j(\Gamma')\cap\tran_j(\mc F_j)$.
We identify $\tran_j(\mc F_j)$ with
$\ab_j(\mc{F})$ as usual.\footnote{An element of $\tran_j(\mc F_j)$
acts on each fiber of
$\mc F_j\to\mc F_{j-1}$ by addition of a fixed
element of $\ab_j(\mc F)$.} 
Since $\Gamma'$ is pure, we have $K_j=\eta_j(\Gamma'_j)$ by \eqref{eq:layer-periods}.
Thus, by \eqref{eq:pure-filtration}, we must compute the
elements added\footnote{By ``elements added'' we mean the following. For the free nilspace $\mc{F}_j$, the induced translation acts by adding one fixed element of the $j$-th structure group. Denote this element by $\delta_j(g)$, that is $\pi_j(g\mf f)
=
\pi_j(\mf f)+\delta_j(g)$ for every $\mf f\in\mc{F}$. 
 The element $\delta_j(g)$ is independent of $\mf f$. These $\delta_j(T_wA_h)$ are the added elements in question.} by $T_wA_h$ when $w\in V_j$ and
$h\in\Gamma_j$.

First we specify the structure groups in which these
elements lie. Write
$F=\prod_{j<k}\mc D_j(\ab_j)$ and
$L_j=\eta_j(\Gamma_j)\le\ab_j$.
The nilspace isomorphism
\eqref{eq:profile-chart} identifies the $j$-th structure
group of $\mc F$, for $j<k$, with
$\ab_j\times(U_j/U_{j+1})$.
Using \eqref{eq:isou} for the second factor gives
\[
 \ab_j(\mc F)\cong
 \ab_j\times(V_j/V_{j+1})\quad \text{ for }j<k.
\]
Both $F$ and $U$ have step at most $k-1$, so only
$\mc D_k(\wt{\ab})$ contributes to the top structure
group. Hence $\ab_k(\mc F)\cong\wt{\ab}$.

Fix $j<k$. We first describe the pair added by $A_h$
for $h\in\Gamma_j$.
Let $\ell_j(h)\in\ab_j$ be the element identified with the constant shift $\eta_j(h)$ (this is indeed constant because $h$ has height at least $j$; see
\cite[Theorem 3.15]{CGSS-doucos}).
By definition of $L_j$, the elements $\ell_j(h)$,
as $h$ ranges over $\Gamma_j$, range over all of $L_j$.

By Lemma~\ref{lem:profile-heights}, $A_h$ also has height
at least $j$. Its induced action on $\mc F_j$
therefore adds a fixed pair
$(\ell_j(h),r_j(h))$, where
$r_j(h)\in V_j/V_{j+1}$ the first component is $\ell_j(h)$ because
$A_h(t,f)=(ht,R_hf)$; indeed applying $A_h$ changes the
$F$-coordinate from $t$ to $ht$.
Thus the element added to the $j$-th coordinate of $F$
is exactly the element $\ell_j(h)$ added by $h$).

Next, we describe how $T_w$, with $w\in V_j$, changes
this pair. Formula \eqref{eq:profile-inverse} shows that
$T_w$ leaves $t$ unchanged and adds $w-w(1)$ to the
$U$-coordinate. Under \eqref{eq:isou}, the resulting
change in the $j$-th structure group of $U$ is
$w+V_{j+1}$: the difference between $w-w(1)$ and $w$
is constant and belongs to $V_{j+1}$ because $j<k$.
Thus $T_w$ adds the pair $(0,w+V_{j+1})$.

We can now determine $K_j$. Every element of $\Gamma'_j$
has the form $T_wA_h$ with $w\in V_j$ and
$h\in\Gamma_j$, so the pair it adds is
$(\ell_j(h),r_j(h)+w+V_{j+1})$.
This proves $K_j\subseteq L_j\times(V_j/V_{j+1})$.
Conversely, take any
$(\ell,r)\in L_j\times(V_j/V_{j+1})$.
Choose $h\in\Gamma_j$ with $\ell_j(h)=\ell$, and then
choose $w\in V_j$ satisfying
$w+V_{j+1}=r-r_j(h)$.
The transformation $T_wA_h$ adds exactly $(\ell,r)$.
This proves the reverse inclusion.

For $j=k$, equation \eqref{eq:pure-filtration} gives
$h\in\Gamma_k=\{1\}$ and $w\in V_k=\Lambda$.
These $w$ are constant functions, and $T_w$ adds their
values to the $\wt{\ab}$ coordinate.
We have therefore proved
\begin{equation}\label{eq:new-periods}
 K_j=L_j\times(V_j/V_{j+1})\quad(j<k),
 \qquad K_k=\Lambda.
\end{equation}

Finally, these formulas verify both required topological
properties. For $j<k$, the subgroup $L_j$ is discrete
by the inductive hypothesis, and $V_j/V_{j+1}$ is
discrete, so $K_j$ is discrete.
Moreover,
$\ab_j(\mc F)/K_j\cong\ab_j/L_j$,
which is compact by the inductive hypothesis.
At the top level, $\Lambda$ is discrete and
$\wt{\ab}/\Lambda\cong\ab$ is compact by
\eqref{eq:top-cover}.
Thus, every $K_j$ is discrete and cocompact, proving
that $\Gamma'$ is fiber-discrete and fiber-cocompact.
\end{proof}

\subsection{Identifying the quotient and completing the pure presentation}\label{subsec:endpurepres}\hfill\\
\begin{proposition}\label{prop:identifying-mcF/Gamma'}
We have the isomorphism of compact nilspaces
\begin{equation}
\mc{F}/\Gamma'\cong \ns.
\end{equation}
\end{proposition}

\begin{proof}
Reduction of functions modulo $\Lambda$ gives the map
\begin{equation}\label{eq:profile-quotient}
\psi_\Lambda:\mc{F}\longrightarrow\mc{E},
 \quad \psi_\Lambda(t,f)=(t,p\co f).
\end{equation}
By \eqref{eq:profile-reduction} we have $p\co f=q_{t,p(f(1))}$.
Thus $\psi_\Lambda(t,f)=(t,q_{t,p(f(1))})\in\mc E$, which verifies that $\psi_\Lambda$ indeed takes values in $\mc{E}\cong\nss$.

Under the identifications \eqref{eq:profile-chart} and \eqref{eq:original-pairs}, we can regard $\psi_\Lambda$ as the map
\[
F\times\mc{D}_k(\wt{\ab})\times U\,
 \to \, F\times\mc{D}_k(\ab),\quad
(t,\wt{z},u)\mapsto(t,p(\wt{z})).
\]
This is an open continuous fibration. It is clearly continuous since $p$ is. To see that it is a fibration, note that $p$ is a surjective homomorphism between abelian groups and thus it induces a fibration between the corresponding degree-$k$ structures $\mc{D}_k(\wt{\ab})$ and $\mc{D}_k(\ab)$. Finally, openness follows from \cite[Theorem 2.21]{CGSS-doucos}.

The fibers of $\psi_\Lambda$ are exactly the $V$-orbits. Indeed, two points $(t,\wt{z},u)$ and $(t',\wt{z}',u')$ in a fiber have $t=t'$ and satisfy $\wt{z}'-\wt{z}\in\Lambda$, so their difference in the function coordinate is $(\wt{z}'-\wt{z})+(u'-u)\in V$. Conversely, addition of a $\Lambda$-valued function leaves $p\co f$ unchanged. The action of $\Gamma'$ is equivariant with $\psi_\Lambda$ via the projection $\Gamma'\to\Gamma$ given by $T_wA_h\mapsto h$. Namely,
\[
\psi_\Lambda(A_h(t,f))=(ht,p\co R_hf)=(ht,R_h(p\co f)) = h\cdot \psi_\Lambda(t, f)
\]
where the action of $\Gamma$ on $\mc E\cong \nss$ is given by \eqref{eq:original-reindexing}.

Now we compose $\psi_\Lambda$ with the factor map from $\mc{E}\cong\nss$ to $\ns$, to obtain a continuous fibration $\psi:\mc{F}\to\ns$. Its fibers are the $\Gamma'$-orbits. If $\psi(\mf{f})=\psi(\mf{f}')$, the images $\psi_\Lambda(\mf{f})$ and $\psi_\Lambda(\mf{f}')$ belong to one lower orbit, so some $h\in\Gamma$ makes $\psi_\Lambda(A_h\mf{f})=\psi_\Lambda(\mf{f}')$. We can then choose an element of $V$ that takes $A_h\mf{f}$ to $\mf{f}'$. Conversely, the map $\psi$ is invariant under the action of $\Gamma'$, i.e. $\psi(T_wA_h\mf{f})=\psi(\mf{f})$. Thus
\begin{equation}\label{eq:quotient-stages}
 \mc{F}/\Gamma'\cong(\mc{F}/V)/\Gamma
 \cong\mc E/\Gamma
 \cong(F\times\mc D_k(\ab))/\Gamma\cong\ns.
\end{equation}
Note that these are also topological nilspace isomorphisms. Indeed, by Proposition~\ref{prop:profile-pure} and \cite[Corollary 5.37]{CGSS-doucos} we have that $\mc F/\Gamma'$ is compact. The induced map $\mc F/\Gamma'\to \ns$ is a continuous bijection, hence a homeomorphism. It preserves cubes because $\psi$ does, and its inverse preserves cubes because every cube of $\ns$ lifts through the fibration $\psi$ to a cube of $\mc{F}$.
\end{proof}
We can now complete the proof of the main result of this section.
\begin{proof}[Proof of Theorem~\ref{thm:pure-presentation}] We proceed by induction on the step $k$ of $\ns$. For $k=0$ the result is trivial and for $k=1$, we have that $\ns$ is isomorphic to $\mc{D}_1(\ab_1(\ns))$ where $\ab_1(\ns)\cong B\times \mb{T}^n$, $B$ is a finite abelian group, and $n\in\mb{Z}_{\ge 0}$. The result then follows by writing $B=\mb{Z}^m/\Lambda_1$, $\mb{T}^n=\mb{R}^n/\mb{Z}^n$, and thus $\ns \cong (\mb{Z}^m\times\mb{R}^n)/(\Lambda_1\times \mb{Z}^n)$.

For $k>1$, assume that $\ns_{k-1}\cong F/\Gamma$ as in \S \ref{sec:orig-exp-pairs}. Then we construct $\mc F$ as in \eqref{eq:profile-space} and $\Gamma'$ as in \eqref{def:Gamma'}. By \eqref{eq:profile-inverse}, \eqref{eq:profile-chart}, and Lemma \ref{lem:U-is-free}, we have that $\mc F$ is (isomorphic to) a free nilspace. By Proposition \ref{prop:profile-pure}, the group $\Gamma'\le\tran(\mc F)$ is a pure translation lattice. The result then follows by Proposition \ref{prop:identifying-mcF/Gamma'}.
\end{proof}

\subsection{Embedding in a toral nilspace}\label{subsec:endembed}\hfill\\
Having proved  Theorem \ref{thm:pure-presentation}, we can now complete the proof of our main embedding theorem.

\begin{proof}[Proof of Theorem~\ref{thm:embedding-intro}]
By Theorem~\ref{thm:pure-presentation}, we have that $\ns\cong F/\Gamma$ where $\Gamma$ is a pure translation lattice. By \cite[Proposition 5.47]{CGSS-doucos}, there is a continuous injective morphism from $\ns$ into a toral nilspace. It follows from the proof of \cite[Proposition 5.47]{CGSS-doucos} that this toral nilspace is also $k$-step.
\end{proof}

\begin{remark}\label{rem:induced}
Recall that a \emph{subnilspace} of a nilspace $\nss$ is a nilspace $\ns$ such that $\cu^n(\ns)\subset \cu^n(\nss)$ for every $n\ge 0$ (in particular, as sets, these nilspaces satisfy $\ns\subset \nss$). We say that $\ns$ is additionally an \emph{induced subnilspace} of $\nss$ if for every cube $\q\in\cu^n(\nss)$ with $\q(v)\in \ns$ for all $v\in \{0,1\}^n$, we have $\q\in \cu^n(\ns)$ (i.e.\ $\cu^n(\nss)\cap \ns^{\{0,1\}^n} = \cu^n(\ns)$). 

The embedding notion involved in Theorem \ref{thm:embedding-intro} does not in itself guarantee that the image of the embedding is an \emph{induced} subnilspace of the toral nilspace. However, this stronger property can indeed be obtained, by exploiting useful features of our specific embedding via pure free nilspace presentations. The only additional tool that is needed for this is presented in Appendix \ref{sec:induced}; see Theorem \ref{thm:inherited-cubes}. A simple application of the latter result yields the following stronger version of Theorem \ref{thm:embedding-intro} (we omit the details).
\end{remark}
\begin{theorem}\label{thm:strongembedding}
 For every $k$-step \textsc{cfr} nilspace $\ns$, there is a $k$-step toral nilspace $\nss$ and an injective continuous morphism $\varphi:\ns\to\nss$ such that $\varphi(\ns)$ is an induced subnilspace of $\nss$.
 Equivalently, every \textsc{cfr} $k$-step nilspace is isomorphic to an induced subnilspace of a $k$-step toral nilspace.
\end{theorem}

\subsection{A concrete 3-step example}\label{subsec:3step}\hfill\\
To demonstrate how the proof manipulates a concrete nilspace with a non-pure presentation, let us revisit Example \ref{ex:problematic-3-step}. We have that $k=3$, $\Gamma\cong\mc{D}_1(\mb{Z}^2)$, and $\Lambda=4\mb{Z}$. Following our proof, we need to consider the space $U$ from \eqref{eq:normalized-U}, whose elements can be checked to be
\[
u(n,m)=4\left(an+bm+c\binom n2+dnm+e\binom m2\right).
\]
Then $a,b$ have weight $2$, while $c,d,e$ have weight $1$. Via the map $\chi$ from \eqref{eq:profile-chart}, we have the nilspace isomorphism
\[
\mc{F}\cong\mc{D}_1(\mb{Z}^2)\times \mc{D}_3(\mb{Z}) \times \mc{D}_2(\mb{Z}^2)
\times\mc{D}_1(\mb{Z}^3).
\]
In coordinates $(x,y,z,a,b,c,d,e)$, the two generators $h_1,h_2$ of $\Gamma$ act on $\mc F$ via $\chi$ as
\begin{align*}
A_{h_1}:&(x,y,z,a,b,c,d,e)\mapsto (x+2,y,z+4a,a+c,b+d-x-1,c,d,e) \quad\text{and}\\
A_{h_2}:&(x,y,z,a,b,c,d,e)\mapsto (x,y+2,z+x^2+4b,a+d,b+e,c,d,e).
\end{align*}
To obtain these formulae, we start with $(x,y,z,a,b,c,d,e)$ and apply $\chi$ to get an element in $\mc{F}$, then apply the definition of $A_h$ from \eqref{eq:profile-actions}, and then go back via $\chi^{-1}$.  

On the other hand, the subgroup $\{T_w:w\in V\}$ is generated by unit shifts in $a,b,c,d,e$ and the shift $z\mapsto z+4$, each fixing all other coordinates (this can be seen similarly using $\chi$).

Together, they exhibit exactly how the auxiliary coordinates absorb the defect in pureness. Note that, in particular, the lattices $\eta_i(\Gamma'_i)$ for $i=1,2,3$ are precisely
$\eta_1(\Gamma'_1)=2\mb{Z}^2\times\mb{Z}^3$, $\eta_2(\Gamma'_2)=\mb{Z}^2$, and $\eta_3(\Gamma'_3)=4\mb{Z}$.

\section{From toral embeddings to the inverse theorem}
We are now ready to prove our main result, Theorem \ref{thm:JT}. The argument follows
that of \cite[Theorem 7.16]{CGSS-doucos}, which we include here
with the details needed for compact abelian groups.

Before we turn to the main proof, we need to address a technical aspect related to the application that we are going to make of the general nilspace inverse theorem  \cite[Theorem 1.6]{inverse}. This theorem as stated involves an underlying assumption that the \textsc{cfr} nilspaces, involved in its correlating nilspace polynomials, all come equipped with a \emph{Riemannian} metric. However, when in this section we use the embedding of those nilspaces into filtered nilmanifolds equipped with Riemannian metrics, we cannot guarantee that the pullback of the latter metric is again Riemannian, nor  bi-Lipschitz equivalent
to the original fixed metric (the embedding is only continuous). Thus the required Lipschitz-bound conclusion in the new pullback metric
does not follow directly from \cite[Theorem 1.6]{inverse}. 
To address this, we first prove a simple upgrade of \cite[Theorem 1.6]{inverse} that holds with arbitrary compatible metrics prescribed on the \textsc{cfr} nilspaces (where ``compatible'' means that the metric generates the nilspace's given topology).

Recall that there are countably many isomorphism classes of $k$-step \textsc{cfr} nilspaces; see \cite[Theorem 2.6.1]{Cand:Notes2}. We choose one $k$-step \textsc{cfr} nilspace from each class and we denote by $\mc{N}_k$ the resulting set of \textsc{cfr} nilspaces. Recall also that for a complex-valued function $F$ on a metric space $X$ with metric $d$, we write $\|F\|_{\Lip}$ for $\|F\|_\infty+\sup_{x,y\in X, x\ne y}
\frac{|F(x)-F(y)|}{d(x,y)}$.

\begin{theorem}[Corollary of Theorem 1.6 from \cite{inverse}]\label{thm:refined-inver-metrics}
Let $k\in \mb N$, and fix an arbitrary compatible metric on each \textsc{cfr} $k$-step nilspace in $\mc{N}_k$. Then, for every $\delta>0$, there exist a finite subset $\mc{N}_{k,\delta}\subset \mc{N}_k$, and $\varepsilon,M>0$ such that the following holds. For every compact abelian group $\ab$ and every 1-bounded Borel function $f:\ab\to\mb{C}$ with $\|f\|_{U^{k+1}}\ge \delta$, there exists $\ns\in\mc{N}_{k,\delta}$, a 1-bounded function $F:\ns\to \mb{C}$ with $\Lip(F)\leq M$, and a continuous morphism $\phi:\mc{D}_1(\ab)\to\ns$ such that $|\mb{E}_{x\in \ab} f(x) F(\phi(x))|\ge \varepsilon$.
\end{theorem}
\begin{proof}
Note that \cite[Theorem 1.6]{inverse} already gives the conclusion provided the chosen metrics on $\mc{N}_k$ are \emph{Riemannian}. For the more general compatible metrics allowed here, we proceed as follows. We first replace the arbitrary chosen metrics by Riemannian ones, and apply \cite[Theorem 1.6]{inverse} to obtain a finite collection $\mc{N}_{k,\delta}$ of target \textsc{cfr} nilspaces $\ns$ and uniform Lipschitz and
correlation bounds $M_0$ and $\varepsilon_0>0$. On each $\ns$ with its fixed Riemannian metric, the family of $1$-bounded functions $F$ with $\Lip(F)\leq M_0$ is compact in the
uniform norm by the Arzel\`a--Ascoli theorem. Lipschitz functions for any
compatible metric are uniformly dense in the continuous
functions. Therefore, finitely many functions that are
Lipschitz for the new metric (and which we can assume to be 1-bounded) approximate this entire family
to within $\varepsilon_0/2$. Taking the maximum of their
Lipschitz constants over the finitely many nilspaces $\ns$ gives a uniform Lipschitz bound for the new metrics. Now given the function $f$ in the theorem, letting $F$ be the Lipschitz function relative to the Riemannian metric such that $|\mb{E}_{x\in \ab} f(x) F(\phi(x))| \geq \epsilon_0$ (as given by \cite[Theorem 1.6]{inverse}), note that we can assume that $F$ is 1-bounded,\footnote{Strictly speaking, \cite[Theorem 1.6]{inverse}) tells us only that $F\co\phi$, but we can easily upgrade this to $F$ itself being 1-bounded. Indeed, project $F:\ns\to \mb{C}$ onto the closed unit disk. This does not increase its Lipschitz constant and preserves the nilspace polynomial $F\co\phi$, which is already $1$-bounded.} and therefore we can replace $F$ by one of the aforementioned approximants, leaving $\phi$ unchanged. This decreases the
correlation by at most $\varepsilon_0/2$, and the result follows.
\end{proof}

\begin{proof}[Proof of Theorem~\ref{thm:JT}]
For each nilspace $\ns\in \mc{N}_k$, applying
Theorem~\ref{thm:embedding-intro} and combining it with \cite[Theorem 2.9.17]{Cand:Notes2}, we obtain an embedding $\iota_\ns$ of $\ns$ into a filtered nilmanifold $G/\Gamma$ of degree at most $k$ (with connected filtration groups), viewed as a $k$-step \textsc{cfr} nilspace. Fix this assignment for the rest of the proof. For each such $\ns$, define the metric
\[
 d_\ns(x,x')
 :=d_{G/\Gamma}\bigl(\iota_\ns(x),\iota_\ns(x')\bigr).
\]
Applying Theorem \ref{thm:refined-inver-metrics} with the metrics $d_\ns$ and parameters $k,\delta$, we obtain $M=M(\delta,k)\in\mb N$, $ \varepsilon= \varepsilon(\delta,k)>0$, a finite collection $\mc{N}_{k,\delta}$ of $k$-step \textsc{cfr} nilspaces, a nilspace 
$\ns\in \mc{N}_{k,\delta}$, a continuous
morphism $\phi:\mc D_1(\ab)\to\ns$, and a $1$-bounded
 function $F_1:\ns\to\mb C$ with $\Lip(F_1)\leq M$ such that $
 |\mb{E}_{x\in \ab} f(x) F_1(\phi(x))|\ge\varepsilon$.

Let $G/\Gamma$ be the chosen nilmanifold into which $\ns$ embeds, and note that it belongs to a finite family $\mc{N}_{k,\delta}'$ of nilmanifolds, obtained by our application of Theorem \ref{thm:embedding-intro} to each $\ns\in \mc{N}_{k,\delta}$. Define $F_2:\iota_\ns(\ns)\to\mb C$ by $F_2(y)=F_1(\iota_\ns^{-1}(y))$. By the definition of $d_\ns$, this $1$-bounded function satisfies $\Lip(F_2)\le M$. We extend this function to $G/\Gamma$ by applying
the McShane--Whitney extension formula separately to its real
and imaginary parts. Define
\begin{multline*}
 F_3(x)
 :=\inf_{y\in\iota_\ns(\ns)}
 \bigl(\Re(F_2(y))+\Lip(\Re(F_2))d_{G/\Gamma}(x,y)\bigr)\\
 {}+i\inf_{y\in\iota_\ns(\ns)}
 \bigl(\Im(F_2(y))+\Lip(\Im(F_2))d_{G/\Gamma}(x,y)\bigr).
\end{multline*}
This function $F_3$ extends $F_2$ and satisfies $
 \Lip(F_3)
 \le \sqrt2\,\Lip(F_1)$.
Further composing with the projection onto the closed
unit disk defined as $z\mapsto z$ for $|z|\le1$ and
$z\mapsto z/|z|$ otherwise, we preserve the extension property and do not increase the Lipschitz constant.
Thus, retaining the notation $F_3$, we have $
 \|F_3\|_\infty\le1$, $
 \Lip(F_3)\le\sqrt2 M$, and $
F_3\co\iota_\ns=F_1$.

Now relabel $\mc{N}_{k,\delta}'$ as $\mc{N}_{k,\delta}$ and set $F:=\overline{F_3}$, and we then have $\varepsilon
\le|\mb{E}_{x\in\ab} f(x) F_1(\phi(x))|
=|\mb{E}_{x\in\ab} f(x) \overline{F\co\iota_\ns\co\phi(x)}|$. The map $\iota_\ns\co\phi$ is a continuous nilspace morphism, that we relabel as $\phi$, and $
\|F\|_{\Lip}\le L:=1+\sqrt2 M$. The family $\mc{N}_{k,\delta}$ of nilmanifolds, equipped with their Riemannian metrics, is finite. This completes the proof.
\end{proof}

\appendix

\section{Toral embeddings with inherited cubes}\label{sec:induced}

The purpose of this appendix is to confirm the main claims in Remark \ref{rem:induced}. Recall from that remark the notion of an induced subnilspace.

\begin{theorem}[Toral embeddings with inherited cubes]
\label{thm:inherited-cubes}
Let $F=\prod_{j=1}^k
 \mc{D}_j(\mb{Z}^{a_j}\times\mb{R}^{b_j})$, and let $\Gamma\leq\tran(F)$ be a pure translation lattice.
Let $\ns=F/\Gamma$, with quotient map $\pi:F\to F/\Gamma$ for the action of $\Gamma$ on $F$. Let $F^*=\prod_{j=1}^k\mc{D}_j(\mb{R}^{a_j+b_j})$, 
and let $\widehat\Gamma$ be the group obtained by extending the
elements of $\Gamma$ to $F^*$ as in
\cite[Definition 5.46 and Proposition 5.47]{CGSS-doucos}, and let $\pi^*:F^*\to F^*/\wh{\Gamma}$ be the corresponding quotient map.
Then $\nss=F^*/\wh{\Gamma}$ is a $k$-step toral nilspace, and the well-defined map $i:\ns\to \nss$ sending $\pi(f)$ to $\pi^*(f)$ is a continuous
injective morphism satisfying
\[
 \forall n\ge 0,\forall\q:\{0,1\}^n\to\ns,\quad \q\in\cu^n(\ns)
 \quad\Longleftrightarrow\quad
 i\co \q\in\cu^n(\nss).
\]
Thus $i$ identifies $\ns$ with a closed induced subnilspace of $\nss$.
\end{theorem}

\begin{proof}
By \cite[Proposition 5.47]{CGSS-doucos}, the quotient $\nss$ is toral, and $i$ is a continuous injective morphism. It remains to prove that $i$ reflects cubes. Identifying $F$ with its coordinate image in $F^*$, we have $\pi^*|_F=i\co \pi$. We use two observations.

First, $F$ is an induced subnilspace of $F^*$: $\cu^n(F)=\cu^n(F^*)\cap F^{\{0,1\}^n}$. Indeed, in each degree-$j$ coordinate, the cube conditions are
alternating-sum equations on $(j+1)$-dimensional faces. For values
in $\mb{Z}^{a_j}\times\mb{R}^{b_j}$, these equations hold
in that subgroup if and only if they hold in
$\mb{R}^{a_j+b_j}$.

Second, every element of $\widehat\Gamma$ preserves $F$ (i.e.\  $ \widehat\gamma t\in F$ for every $t\in F$, $ \widehat\gamma\in \widehat\Gamma $) and
restricts there to the corresponding element of $\Gamma$. Therefore $F$ is a union of entire $\widehat\Gamma$-orbits,
and hence
\[
 (\pi^*)^{-1}(i(\ns))=(\pi^*)^{-1}(i(\pi(F)))=(\pi^*)^{-1}(\pi^*(F))=F.
\]
Now suppose that $\q:\cube{n}\to \ns$ satisfies $i\co \q\in\cu^n(\nss)$. By definition of the quotient cube structure on $\nss$, there is
$\widetilde \q\in\cu^n(F^*)$ such that $
 \pi^*\co\wt\q=i\co \q$. 
The second observation forces every vertex of $\wt{\q}$ to lie in $F$. The first observation therefore gives
$\wt{\q}\in\cu^n(F)$. Thus $\pi\co\wt{\q}$ is a cube
of $\ns$, and injectivity of $i$ implies
$\pi\co\widetilde \q=\q$ (since $i\co \pi\co \wt{\q}=\pi^*\co\wt\q=i\co \q$). This proves the required equivalence.

Since $\ns$ is compact and $\nss$ is Hausdorff, $i$ is a homeomorphism
onto a closed image.
\end{proof}

\noindent \textbf{Funding.} All authors used funding from project PID2024-156180NB-I00 funded by the MICIU/AEI and the European Union. The second-named author was funded by HORIZON-MSCA-2024-PF-01, AlgHOF 101202161.\footnote{Funded by the European Union. Views and opinions expressed are those of the author(s) only and do not reflect those of the European Union or the European Commission. Neither the European Union nor the European Commission can be held responsible for them.} The third-named author was supported by the NKFIH excellence 154126 grant.

\medskip

\noindent \textbf{Disclosure on LLM use.} The authors have been working on the Jamneshan--Tao conjecture for several years and have developed various techniques in the process. One of the most important techniques to emerge from this effort is based on their result showing that, in a pure presentation of a nilspace as a factor of a free nilspace, the discrete coordinates can be replaced by continuous ones to obtain an embedding of the nilspace into a nilmanifold \cite[Proposition 5.47]{CGSS-doucos}.

The authors realized that embedding arbitrary \textsc{cfr} nilspaces into nilmanifolds would imply the Jamneshan–Tao conjecture. This approach gained new momentum recently, when the authors proved the conjecture for the $U^4$-norm using a special case of such embeddings (unpublished), and also when they found that large classes of 3-step \textsc{cfr} nilspaces satisfy the conclusion of Theorem \ref{thm:embedding-intro}, which gave renewed hope that this theorem might hold in general.

LLMs entered the research at this stage. The authors provided the GPT-Astra model with their results on the $U^4$-norm and worked with it toward a generalization. This collaborative process produced an initial version of the proof, which was subsequently refined by the authors in conversation with GPT-Astra through many iterations.

\end{document}